\documentclass[12pt,reqno]{amsart}
\usepackage{amssymb}
\usepackage{amscd}
\usepackage{amsmath}

\newfont{\cyr}{wncyr10} 

\theoremstyle{plain}
\newtheorem{theorem}{Theorem}[section]
\newtheorem{Thm}{Theorem}
\newtheorem{lemma}[theorem]{Lemma}
\newtheorem{proposition}[theorem]{Proposition}
\newtheorem{corollary}[theorem]{Corollary}

\newtheorem*{proposition*}{Proposition}

\theoremstyle{definition}
\newtheorem{remark}[theorem]{Remark}
\newtheorem{definition}[theorem]{Definition}

\numberwithin{equation}{section}

\newcommand{\cB}{{\mathcal B}}

\newcommand{\cC}{{\mathcal C}}
\newcommand{\bC}{{\mathbf C}}

\newcommand{\bD}{{\mathbf D}}

\newcommand{\cE}{{\mathcal E}}

\newcommand{\cF}{{\mathcal F}}
\newcommand{\ff}{{\mathfrak f}}
\newcommand{\fF}{{\mathfrak F}}

\newcommand{\cG}{{\mathcal G}}

\newcommand{\ch}[1]{{\check{#1}}}
\newcommand{\cH}{{\mathcal{H}}}

\newcommand{\cI}{{\mathcal I}}
\newcommand{\cJ}{{\mathcal J}}

\newcommand{\1}{{\mathbf 1}}  

\newcommand{\cK}{{\mathcal K}}
\newcommand{\fK}{{\mathfrak K}}
\renewcommand{\L}{{\mathcal L}}
\newcommand{\cL}{{\mathcal L}}

\newcommand{\cM}{{\mathcal M}}

\newcommand{\ms}{{\Sigma}}  
\newcommand{\bN}{{\mathbf N}}
\newcommand{\cN}{{\mathcal N}}
\newcommand{\cO}{{\mathcal O}}
\newcommand{\ov}[1]{{\overline{#1}}}

\newcommand{\fp}{{\mathfrak p}}

\newcommand{\bQ}{{\mathbf Q}}
\newcommand{\fq}{{\mathfrak q}}
\newcommand{\cR}{{\mathcal R}}

\newcommand{\rel}{\mathrm {rel}}

\newcommand{\sha}{{\mbox{\cyr Sh}}}

\newcommand{\str}{\mathrm {str}}
\newcommand{\vt}{{\vartheta}}

\newcommand{\cX}{{\mathcal X}}
\newcommand{\cY}{{\mathcal Y}}
\newcommand{\bZ}{{\mathbf Z}}

\DeclareMathOperator{\Aut}{Aut}
\DeclareMathOperator{\cha}{Char}

\DeclareMathOperator{\Coker}{Coker}

\DeclareMathOperator{\Div}{div}

\DeclareMathOperator{\Gal}{Gal}

\DeclareMathOperator{\Hom}{Hom}
\DeclareMathOperator{\Image}{Im}
\DeclareMathOperator{\Ker}{Ker}

\DeclareMathOperator{\loc}{\mathrm{loc}}

\DeclareMathOperator{\ord}{ord}

\DeclareMathOperator{\rk}{rk}

\DeclareMathOperator{\Sel}{Sel}

\DeclareMathOperator{\tors}{tors}

\DeclareMathOperator{\Tw}{Tw}

\begin{document}
\title[Birch and Swinnerton-Dyer conjecture]{On Rubin's variant
of the $p$-adic Birch and\\ Swinnerton-Dyer conjecture}
\author{A. Agboola}
\date{Final version. January 22, 2007}
\address{Department of Mathematics \\
University of California \\ Santa Barbara, CA 93106. }
\email{agboola@math.ucsb.edu}
\thanks{Partially supported by NSF grant DMS-0401319.}
\begin{abstract}
We study Rubin's variant of the $p$-adic Birch and Swinnerton-Dyer
conjecture for CM elliptic curves concerning certain special values of
the Katz two-variable $p$-adic $L$-function that lie outside the range
of $p$-adic interpolation.
\end{abstract}

\maketitle

\section{Introduction} \label{S:intro}

Let $E/\bQ$ be an elliptic curve with complex multiplication by $O_K$,
the ring of integers of an imaginary quadratic field $K$ (necessarily
of class number one). Let $p>3$ be a prime of good, ordinary reduction
for $E$; then we may write $pO_K = \fp \fp^*$, with $\fp = \pi O_K$
and $\fp^* = \pi^* O_K$.

Set $\cK_\infty:= K(E_{\pi^{\infty}})$, $\cK^*_\infty:=
K(E_{\pi^{*\infty}})$, and $\fK_\infty:= \cK_\infty
\cK^*_\infty$. Write $K_\infty$ (resp. $K^*_\infty$) for the unique
$\bZ_p$ extension of $K$ unramified outside $\fp$ (resp. $\fp^*$). Let
$\cO$ denote the completion of the ring of integers of the maximal
unramified extension of $\bQ_p$. For any extension $L/K$ we set
$\Lambda(L):= \Lambda(\Gal(L/K)):= \bZ_p[[\Gal(L/K)]]$, and
$\Lambda(L)_{\cO}:= \cO[[\Gal(L/K)]]$.  We write $X(L)$
(resp. $X^*(L)$) for the Pontryagin dual of the $\fp$-primary Selmer
group $\Sel(L, E_{\pi^\infty})$ (resp. the $\fp^*$-primary Selmer
group $\Sel(L, E_{\pi^{*\infty}})$) of $E/L$.

Let
\begin{align*}
&\psi: \Gal(\ov{K}/K) \to \Aut(E_{\pi^{\infty}}) \xrightarrow{\sim}
O_{K,\fp}^{\times} \xrightarrow{\sim} \bZ_{p}^{\times}, \\
&\psi^*: \Gal(\ov{K}/K) \to \Aut(E_{\pi^{*\infty}}) \xrightarrow{\sim}
O_{K,\fp^*}^{\times} \xrightarrow{\sim} \bZ_{p}^{\times}
\end{align*}
denote the natural $\bZ_{p}^{\times}$-valued characters of
$\Gal(\ov{K}/K)$ arising via Galois action on $E_{\pi^\infty}$ and
$E_{\pi^{*\infty}}$ respectively. We may identify $\psi$ with the
Grossecharacter associated to $E$ (and $\psi^*$ with the complex
conjugate $\ov{\psi}$ of this Grossencharacter), as described, for
example, in \cite[p. 325]{R}. We write $T$ (resp. $T^*$) for the
$\fp$-adic (resp. $\fp^*$-adic) Tate module of $E$.

The two-variable Iwasawa main conjecture (proved by Rubin \cite{R2})
implies that $X(\fK_{\infty})$ is a torsion
$\Lambda(\fK_{\infty})$-module whose characteristic ideal in
$\Lambda(\fK_{\infty})_{\cO}$ is generated by a twist of Katz's
two-variable $p$-adic $L$-function $\L_{\fp}$ by the character
$\psi$. The function $\L_{\fp}$ satisfies a $p$-adic interpolation
formula that may be described as follows (see \cite[Theorem 7.1]{R}
for the version given here, and also \cite[Theorem II.4.14]{dS}). 
For all pairs of integers $j,k \in \bZ$ with $0 \leq -j <k$, and for
all characters $\chi: \Gal(K(E_p)/K) \to \ov{K}^{\times}$, we have
\begin{equation} \label{E:interpolate}
\L_{\fp}(\psi^k \psi^{*j} \chi) = A \cdot L(\psi^{-k} \ov{\psi}^{-j}
\chi^{-1}, 0).
\end{equation}
Here $L(\psi^{-k} \ov{\psi}^{-j} \chi^{-1},s)$ denotes the complex
Hecke $L$-function, and $A$ denotes an explicit, non-zero factor whose
precise description need not concern us here.

Define
$$
L_{\fp}(s):= \L_{\fp}(\psi <\psi>^{s-1}),\quad L_{\fp}^{*}(s):=
\L_{\fp}(\psi^* <\psi^{*}>^{s-1})
$$
for $s \in \bZ_p$. The character $\psi$ lies within the range of
interpolation of $\cL_{\fp}$, and the $\fp$-adic Birch and
Swinnerton-Dyer conjecture for $E$ (see \cite[pages 133--134]{BGS},
\cite[Theorem V.8]{PR7}) predicts that $\ord_{s=1} L_{\fp}(s)$ is
equal to the rank $r$ of $E(\bQ)$, and that
\begin{equation*} 
\lim_{s \to 1} \frac{L_{\fp}(s)}{(s-1)^r} \sim
[\log_p(\psi(\gamma_1))]^r \cdot
\left(1 - \frac{\psi(\fp)}{p} \right)\cdot 
\left(1 - \frac{\psi(\fp^*)}{p} \right)\cdot 
|\sha(K)(\fp)| \cdot R_{K,\fp},
\end{equation*}
where $\gamma_1$ is a topological generator of $\Gal(\cK_\infty/K)$,
$\sha(K)(\fp)$ is the $\fp$-primary component of the Tate-Shafarevich
group $\sha(K)$ of $E/K$, $R_{K,\fp}$ is the regulator associated to
the algebraic $\fp$-adic height pairing
$$
\{\,, \,\}_{K,\fp}: \Sel(K,T^*) \times \Sel(K,T) \to
O_{K,\fp}
$$
on $E/K$ (see \cite{PR1}), and the symbol `$\sim$' denotes equality up
to multiplication by a $p$-adic unit.

On the other hand, the character $\psi^*$ lies outside the range of
interpolation of $\L_{\fp}$ and the function $L_{\fp}^*(s)$ has not
been studied nearly as much as $L_{\fp}(s)$. The only results
concerning $L_\fp^*(s)$ of which the author is aware are due to Rubin
(see \cite{R}, \cite{R1}). When $r \geq 1$, Rubin formulated a variant
of the $\fp$-adic Birch and Swinnerton-Dyer conjecture for
$L_\fp^*(s)$ which predicts that that $\ord_{s=1} L_{\fp}^*(s)$ is
equal to $r-1$, and which gives a formula for $\lim_{s \to
1}[L_{\fp}^*(s)/(s-1)^{r-1}]$. Under suitable hypotheses, Rubin showed
that his conjecture is equivalent to the usual $\fp$-adic Birch and
Swinnerton-Dyer conjecture, and he proved both conjectures when
$r=1$. In the case $r=1$, he then used these results to give a
striking $p$-adic construction of a global point of infinite order in
$E(\bQ)$ directly from the special value of a $p$-adic $L$-function.

When $r=0$, however, the above analysis breaks down, and the situation
is less clear. The functional equation satisfied by $\cL_\fp$ (see
\cite[II \S6]{dS}) shows that $\ord_{s=1} L_\fp(s)$ and $\ord_{s=1}
L_\fp^*(s)$ have opposite parity, and so when $r=0$, one expects
that $\ord_{s=1} L_\fp^*(s)$ is odd. This may perhaps be viewed as
being an analogue of a similar exceptional zero phenomenon observed in
the work of Mazur, Tate and Teitelbaum concerning $p$-adic Birch and
Swinnerton-Dyer conjectures for elliptic curves \textit{without}
complex multiplication (see \cite{MTT}, \cite{Gr94}). As Rubin points
out (see \cite[Remark on p. 74]{R1}), it is reasonable to guess that
$\ord_{s=1} L_\fp^*(s)=1$. If this is so, then one would like to
determine the value of $\lim_{s \to 1} [L_\fp^*(s)/(s-1)]$.

In this paper we study an Iwasawa module naturally associated to
$L_\fp^*(s)$ via the two-variable main conjecture and, among other
things, we prove that the above guess is indeed correct. The Iwasawa
module in question is the Pontryagin dual $X_{\fp^*}(K_\infty^*,W^*)$
of a certain \textit{restricted Selmer group}
$\ms_{\fp^*}(K_\infty^*,W^*)$. This restricted Selmer group is defined
by \textit{reversing} the Selmer conditions above $\fp$ and $\fp^*$
that are used to define the usual Selmer group
$\Sel(K_\infty^*,W^*)$. The two-variable main conjecture implies that
a characteristic power series $H_K \in \Lambda(K_\infty^*)$ of
$X_{\fp^*}(K_\infty^*,W^*)$ may be viewed as being an algebraic
$p$-adic $L$-function corresponding to $L_{\fp}^{*}(s)$. We study
$L_{\fp}^{*}(s)$ by analysing the behaviour of $H_K$.

A special case of our results may be described as follows. We define a
compact restricted Selmer group $\ch{\ms}_{\fp^*}(K,T^*) \subseteq
H^1(K,T^*)$. The $O_{K,\fp^*}$-module $\ch{\ms}_{\fp^*}(K,T^*)$ is
free of rank $|r-1|$, and if $r \geq 1$, then it lies in the usual
Selmer group $\Sel(K,T^*)$ associated to $T^*$. The $O_{K,\fp^*}$-rank
of $\ch{\ms}_{\fp^*}(K,T^*)$ governs the order of vanishing of
$L_{\fp}^{*}(s)$ at $s=1$ in the same way that the $O_{K,\fp}$-rank of
$\Sel(K,T)$ determines $\ord_{s=1} L_{\fp}(s)$. We also define a
similar group $\ch{\ms}_\fp(K,T) \subseteq H^1(K,T)$, and we explain
how to construct a $p$-adic height pairing
$$
[\,, \,]_{K,\fp^*}: \ch{\ms}_\fp(K,T) \times \ch{\ms}_{\fp^*}(K,T^*) \to
O_{K,\fp^*}.
$$ 
If $r \geq 1$, then in fact $\ch{\ms}_\fp(K,T) \subseteq \Sel(K,T)$,
$\ch{\ms}_{\fp^*}(K,T^*) \subseteq \Sel(K,T^*)$, and, if the
$\fp^*$-adic Birch and Swinnerton-Dyer conjecture is true, then the
$p$-adic height pairing $[\,, \,]_{K,\fp^*}$ is non-degenerate. We
conjecture that $[\,, \,]_{K,\fp^*}$ is also non-degenerate when $r=0$
(see Remark \ref{R:height}).

Define
$$
\sha_{\rel(\fp)}(K):= \Ker \left[ H^1(K,E) \to \prod_{v \nmid \fp}
H^1(K_v,E) \right],
$$ 
and write $\sha_{\rel(\fp)}(K)(\fp^*)$ for its $\fp^*$-primary
subgroup. Let $\sha_{\rel(\fp)}(K)(\fp^*)_{/\Div}$ denote the quotient
of $\sha_{\rel(\fp)}(K)(\fp^*)$ by its maximal divisible subgroup. It
may be shown that $\sha_{\rel(\fp)}(K)(\fp^*)$ has
$O_{K,\fp^*}$-corank one, and that
$\sha_{\rel(\fp)}(K)(\fp^*)_{/\Div}$ is finite.

\begin{Thm} \label{T:A}
Suppose that $[\, , \,]_{K,\fp^*}$ is non-degenerate, and let
$\gamma$ be a topological generator of $\Gal(\cK_\infty^*/K)$. Then,
if $r=0$, we have $\ord_{s=1} L_\fp^*(s)=1$, and
\begin{align*} 
\lim_{s \to 1}& \frac{L_\fp^*(s)}{s-1} \sim \\
&\log_p(\psi^*(\gamma)) \cdot
(1 - \psi(\fp^*)) \cdot
\frac{|\sha_{\rel(\fp)}(K)(\fp^*)_{/\Div}|}{[H^1(K_{\fp^*},T):
\loc_{\fp^*}(\ms_{\fp}(K,T)]} \cdot \cR_{K,\fp^*},
\end{align*}
where $\cR_{K,\fp^*}$ is a $p$-adic regulator associated to $[\, ,
\,]_{K,\fp^*}$.
\end{Thm}

We also obtain an exact (but much less explicit) formula for $\lim_{s
  \to 1} L^{*}_{\fp}(s)/(s-1)$ by applying the methods of \cite{R} in
  our present setting (see Theorem \ref{T:exact} below).

Suppose now that $r \geq 1$, and assume that $\sha(K)(p)$ is
finite. Then $E(K) \otimes_{O_K} O_{K,\fp^*}$ is a free
$O_{K,\fp^*}$-module of rank $r$, and the kernel of the localisation
map
$$
E(K) \otimes_{O_K} O_{K,\fp^*} \to E(K_{\fp^*}) \otimes_{O_K} O_{K,\fp^*}
$$
has $O_{K,\fp^*}$-rank $r-1$. Let $y_1, \ldots, y_{r-1}$ be an
$O_{K,\fp^*}$-basis of this kernel, and extend it to an
$O_{K,\fp^*}$-basis $y_1, \ldots, y_{r-1},y_{\fp^*}$ of $E(K)
\otimes_{O_K} O_{K,\fp^*}$. We write $x_1,\ldots, x_{r-1}, y_{\fp}$
for a similarly constructed $O_{K,\fp}$-basis of $E(K) \otimes_{O_K}
O_{K,\fp}$. The following result is a direct consequence of Rubin's
precise formula for $\lim_{s \to 1} [L_{\fp}^{*}(s)/(s-1)^{r-1}]$ (see
\cite[Corollary 11.3]{R}). We give a new proof of this result which is
different from that contained in \cite{R}. In particular, our proof
gives an alternative way of viewing the somewhat unusual regulator
$R_{\fp}^{*}$ defined in \cite[\S11]{R}.

\begin{Thm} \label{T:B}
Suppose that $r \geq 1$ and that $[\,, \,]_{K,\fp^*}$ is
non-degenerate. Then $\ord_{s=1} L_\fp^*(s) = r-1$, and 
\begin{align} \label{E:bsd3}
\lim_{s \to 1}& \frac{L_\fp^*(s)}{(s-1)^{r-1}} \sim \notag \\
&[\log_p(\psi^*(\gamma))]^{r-1} \cdot p^{-2} \cdot
|\sha(K)(\fp^*)| \cdot
\log_{E,\fp^*}(y_{\fp^*}) \cdot \log_{E,\fp}(y_{\fp})
\cdot \cR_{K,\fp^*},
\end{align}
where $\log_{E,\fp^*}$ (resp. $\log_{E,\fp}$) denotes the $\fp^*$-adic
(resp. $\fp$-adic) logarithm associated to $E$.
\end{Thm}

An outline of the contents of this paper is as follows. In Section
\ref{S:twists} we recall some basic facts about twists of Iwasawa
modules and derivatives of characteristic power series, and we apply
these results to describe the relationship between $L_{\fp}^{*}(s)$
and a characteristic power series $H_K \in \Lambda(K_\infty^*)$ of
$X_{\fp^*}(K_\infty^*,W^*)$. In Section \ref{S:selmer} we define
various Selmer groups, and we establish some of their properties. We
describe how to construct an algebraic $p$-adic height pairing on
restricted Selmer groups in Section \ref{S:height}. In Section
\ref{S:leading} we calculate (under certain hypotheses) the leading
term of a characteristic power series $H_F \in \Lambda(F_\infty^*)$ of
$X_{\fp^*}(F_\infty^*,W^*)$, where $F/K$ is any finite extension, and
$F_\infty^*:= FK_{\infty}^{*}$. In Section \ref{S:rsk} we study
restricted Selmer groups over $K$, and we show that, under certain
standard assumptions, $\ord_{s=1} L_{\fp}^{*}(s) = |r-1|$. We then
give the proof of Theorem \ref{T:A} in Section \ref{S:ThmA}, and that
of Theorem \ref{T:B} in Section \ref{S:ThmB}. Finally, in Section
\ref{S:exact}, we explain how the methods of \cite{R} may be used to
give a formula for the exact value of $\lim_{s \to 1}
L_{p}^{*}(s)/(s-1)$ when $r=0$.
\medskip

\noindent{}{\bf Acknowledgements.} I am very grateful indeed to Karl
Rubin for extremely helpful conversations and correspondence. Parts of
this paper were written while I was visiting the Universit\'e de
Bordeaux I and the Centre de Recherches Math\'ematiques at the
Universit\'e de Montreal. I thank these institutions for their
hospitality and support.
\medskip

\noindent{}{\bf Notation and conventions.}
For each integer $n \geq 1$, we write
$$
\cK_n:= K(E_{\pi^n}),\quad \cK_n^*:= K(E_{\pi^{*n}}).
$$ 

For each place $v$ of $K$, we write $k_v$ for the residue field of
$v$, and $\tilde{E}_v/k_v$ for the reduction of the elliptic curve $E$
modulo $v$. We set $W:= E_{\pi^\infty}$ and $W^*:= E_{\pi^{*\infty}}$.

Throughout this paper, $F$ denotes a finite extension of $K$, and we
set
\begin{align*}
&\cF_n:= F\cK_n,\quad \cF_\infty:= F\cK_\infty,\quad F_\infty:=
FK_\infty,\\
&\cF_n^*:= F\cK_n^*,\quad \cF_\infty^*:= F\cK_\infty^*,\quad F_\infty^*:=
FK_\infty^*, \\
&\fF_\infty:= F\fK_\infty.
\end{align*}

For any extension $L/K$ we write $\cM(L)$ (resp. $\cM^*(L)$) for the
maximal abelian pro-$p$ extension of $L$ which is unramified away from
$\fp$ (resp. $\fp^*$), and we set
$$
\cX(L):= \Gal(\cM(L)/L),\quad \cX^*(L):= \Gal(\cM^*(L)/L).
$$
We let $\cB(L)$ (resp. $\cB^*(L)$) denote the maximal abelian pro-$p$
extension of $L$ which is unramified away from $\fp$ (resp. $\fp^*$)
and totally split at all places of $L$ lying above $\fp^*$
(resp. $\fp$), and we write
$$
\cY(L):= \Gal(\cB(L)/L),\quad \cY^*(L):= \Gal(\cB^*(L)/L).
$$

If $M$ is any $\bZ_p$-module, then $M_{\Div}$ denotes the maximal
divisible submodule of $M$, and we set $M_{/\Div}:= M/M_{\Div}$. We
write $M_{\tors}$ for the torsion submodule of $M$, and $M^\land$ for
the Pontryagin dual of $M$. If $M$ is a torsion $O_{K,\fq}$-module,
with $\fq \in \{ \fp, \fp^* \}$, then we write $T_{\fq}(M)$ for the
$\fq$-adic Tate module of $M$.

We set $D_\fp:= K_\fp/O_{K,\fp}$ and $D_{\fp^*}:= K_{\fp^*}/O_{K,\fp^*}$.

\section{Twists and derivatives} \label{S:twists}

In this section we shall recall some basic facts concerning twists of
Iwasawa modules and derivatives of characteristic power series. We
then apply these results to a twist of the Katz two-variable
$p$-adic $L$-function $\cL_\fp$ by the character $\psi^*$.

Let $\cG_F:= \Gal(\fF_\infty/F)$, and suppose that $\rho: \cG_F \to
\bZ_p^\times$ is any character. Then we have a twisting map
$$
\Tw_\rho: \Lambda(\cG_F) \to \Lambda(\cG_F)
$$
associated to $\rho$ which is induced by the map $g \mapsto \rho(g)g$
for all $g \in \cG_F$. If $M$ is a finitely generated
$\Lambda(\cG_F)$-module with characteristic power series $f_M$, then a
routine computation shows that $\Tw_\rho(f_M)$ is a characteristic
power series of $M(\rho^{-1}):= M \otimes \rho^{-1}$.

Set $\cH:= \Ker(\rho)$. Then there is a natural quotient map
$$
\Pi_{\cG_F/\cH}: \Lambda(\cG_F) \to \Lambda(\cG_F/\cH),
$$
and $\Pi_{\cG_F/\cH}(\Tw_\rho(f_M))$ is a characteristic power series of
the $\Lambda(\cG_F/\cH)$-module $M(\rho^{-1}) \otimes_{\Lambda(\cG_F)}
\Lambda(\cG_F/\cH)$. If $\rho_1: \cG_F \to \bZ_p^\times$ is any character
which factors through $\cG_F/\cH$, then
\begin{equation} \label{E:LT1a}
[\Tw_{\rho}(f_M)](\rho_1) = [\Pi_{\cG_F/\cH}(\Tw_\rho(f_M))](\rho_1),
\end{equation}
and there is an isomorphism
$$
M(\rho^{-1}) \otimes_{\Lambda(\cG_F)} \Lambda(\cG_F/\cH) \simeq
(M \otimes_{\Lambda(\cG_F)} \Lambda(\cG_F/\cH))(\rho^{-1})
$$
of $\Lambda(\cG_F/\cH)$-modules. Hence we may study the values of
$\Tw_\rho(f_M)$ at characters $\rho_1$ which factor through $\cG_F/\cH$
by studying the values of $\Pi_{\cG/\cH}(\Tw_\rho(f_M))$ at such
characters.

Suppose now that $\rho$ is of infinite order, and let $N$ be a
finitely generated $\Lambda(\cG_F/\cH)$-module with characteristic power
series $f_N \in \Lambda(\cG_F/\cH)$. We may write
$$
\cG_F/\cH \simeq \Delta \times G,
$$
where $|\Delta|$ is prime to $p$, and $G \simeq \bZ_p$. Let $\gamma$
be a fixed topological generator of $\cG_F/\cH$, and let $\Pi_G:
\Lambda(\cG_F/\cH) \to \Lambda(G)$ be the natural quotient map. We
identify $\Lambda(G)$ with $\bZ_p[[t]]$ in the usual way via the map
$\Pi_G(\gamma) \mapsto 1+t$.

Let $I_{\cG_F/\cH}$ denote the augmentation ideal of $\Lambda(\cG_F/\cH)$,
and suppose that $n \geq 0$ is the largest integer such that $f_N \in
I_{\cG_F/\cH}^{n}$ and $f_N \notin I_{\cG_F/\cH}^{n+1}$. It is not hard to
check that $\Pi_G(f_N)(t)$ is a characteristic power series of the
$\Lambda(G)$-module $N^{\Delta}$, and that
\begin{equation} \label{E:LT1}
((\gamma -1)^{-n}f_N)(\1) = \frac{\Pi_G(f_N)}{t^n} \Biggr |_{t=0},
\end{equation}
where $\1$ denotes the identity character of $\cG_F/\cH$.

For any character $\nu: \cG_F/\cH \to \bZ_p^\times$, we set
$\vartheta_\nu:= \nu(\gamma)^{-1} \gamma - 1$. Then if $m \geq 0$ is
any integer, it follows from the definitions that we have
\begin{equation} \label{E:LT2}
(\vartheta^{-m}_{\nu} f_N)(\nu) = [(\gamma-1)^{-m} \Tw_\nu(f_N)] (\1),
\end{equation}
where $\Tw_\nu: \Lambda(\cG_F/\cH) \to \Lambda(\cG_F/\cH)$ is the twisting
map associated to $\nu$.

We now recall how \eqref{E:LT2} is related to derivatives
of certain $p$-adic analytic functions as described in
\cite[\S7]{R}. Write $<\nu>: \cG_F/\cH \to \bZ_p^\times$ for the
composition of $\nu$ with the natural projection $\bZ_p^\times \to 1
+p\bZ_p$, and suppose that $\chi: \cG_F/\cH \to \bZ_p^\times$ is any
character of order prime to $p$. The map from $\bZ_p$ to $\bC_p$
given by $s \mapsto f_N(\nu \chi <\nu>^{s-1})$ defines an analytic
function on $\bZ_p$. Define
$$
\ord_{\nu \chi}(f_N):= \ord_{s=1} f_N(\nu \chi <\nu>^{s-1}),
$$
and set
$$
\bD^{(m)}f_N(\nu \chi):= \frac{1}{m!} \left( \frac{d}{ds} \right)^m f_N(\nu
\chi <\nu>^{s-1}) \Biggr |_{s=1}.
$$
We write
$$
f_{N}^{(m)}(\nu \chi):= \bD^{(m)} f_{N}(\nu \chi),
$$
and we extend these definitions to $\Lambda(\cG_F)$ via the quotient map
$\Pi_{\cG_F/\cH}$. A routine calculation shows that we have
$$
\bD^{(m)}(\vartheta_\nu^m(\nu \chi)) = \{ \log_p(\nu(\gamma)) \}^m,
$$
and
\begin{equation} \label{E:LT3}
\bD^{(m)}(\vartheta_\nu^m f_N)(\nu \chi) =
\{ \log_p(\nu(\gamma)) \}^m f_N(\nu \chi) =
[\{ \log_p(\nu(\gamma)) \}^m \Tw_\nu(f_N)](\chi).
\end{equation}

We can now see from \eqref{E:LT1}, \eqref{E:LT2} and \eqref{E:LT3}
that if $n_\nu:= \ord_\nu(f_N)$, then we may write $f_N =
\vartheta_{\nu}^{n_{\nu}} F_\nu$ with $F_\nu \in \Lambda(\cG_F/\cH)$,
and we have
\begin{align} \label{E:LT4}
f_{N}^{(n_{\nu})}(\nu) &= \lim_{s \to 1} \frac{f_N(\nu
<\nu>^{s-1})}{(s-1)^{n_\nu}} \notag \\
&= \bD^{(n_\nu)} ( \vartheta_{\nu}^{n_{\nu}} F_\nu
) (\nu) \notag \\
&= [ \{ \log_p(\nu(\gamma)) \}^{n_\nu} \Tw_\nu(F_\nu)] (\1) \notag \\
&= \{ \log_p(\nu(\gamma)) \}^{n_\nu} \cdot \Pi_G(\Tw_\nu(F_\nu))(0)
\notag \\
&= \{ \log_p(\nu(\gamma)) \}^{n_\nu} \cdot
\frac{\Pi_G(\Tw_\nu(f_N))}{t^{n_\nu}} \Biggr |_{t=0}.
\end{align}

We shall now apply the above discussion to the case in which $F=K$,
$M=\cX(\fK_\infty)$, $\rho = \nu = \psi^*$, $\cH=
\Gal(\fK_\infty/\cK^*_\infty)$, $G = \Gal(K_\infty^*/K)$ and $\chi =
\1$.

Recall that the two-variable main conjecture asserts that
$\cX(\fK_\infty)$ is a torsion $\Lambda(\fK_\infty)$-module, and that
the Katz two-variable $p$-adic $L$-function $\cL_\fp$ is a
characteristic power series of $\cX(\fK_\infty)$ in
$\Lambda(\fK_\infty)_{\cO}$. We therefore see that
$\Tw_{\psi^*}(\L_{\fp}) \in \Lambda(\fK_\infty)_{\cO}$ is a
characteristic power series of $\cX(\fK_\infty)(\psi^{*-1})$. Let
$I_{K_\infty^*}$ denote the kernel of the natural map
$\Lambda(\fK_\infty) \to \Lambda(K_\infty^*)$. Fix any characteristic
power series $H_K \in \Lambda(K^*_\infty)$ of the
$\Lambda(K^*_\infty)$-module
$$
\cX(\fK_\infty)(\psi^{*-1})
\otimes_{\Lambda(\fK_\infty)}
(\Lambda(\fK_\infty)/I_{K^*_\infty}) \simeq
\cX(\fK_\infty)(\psi^{*-1})/I_{K^*_\infty}\cX(\fK_\infty)(\psi^{*-1}).
$$
Then we deduce from \eqref{E:LT1a}, \eqref{E:LT1} and \eqref{E:LT4}
that
\begin{equation} \label{E:orders}
\ord_{s=1} L_\fp^*(s) = \ord_{t=0} H_K,
\end{equation}
and if we set $n_{\psi^*}:= \ord_{s=1} L_\fp^*(s)$, then
\begin{equation} \label{E:LT6}
\cL_{\fp}^{(n_{\psi^*})} (\psi^*) = \lim_{s \to 1}
\frac{L_\fp^*(s)}{(s-1)^{n_{\psi^*}}} \sim
\{ \log_p(\psi^*(\gamma)) \}^{n_{\psi^*}} \cdot
\frac{H_K}{t^{n_{\psi^*}}} \Biggr |_{t=0},
\end{equation}
where `$\sim$' denotes equality up to multiplication by a $p$-adic
 unit (in fact, in this case, we have equality up to multiplication by
 an element of $\cO^\times$).

\section{Selmer Groups} \label{S:selmer}

In this section we shall define various Selmer groups that we require,
and establish some of their properties.

For any place $v$ of $F$, we define
$H^1_f(F_v,W)$ to be the image of $E(F_v) \otimes D_\fp$ under
the Kummer map
$$
E(F_v) \otimes D_\fp \to H^1(F_v,W),
$$
and we define $H^1_f(F_v,W^*)$ in a similar manner. Note that
$H^1_f(F_v,W)=0$ if $v\nmid \fp$. We also set
\begin{align*}
&H^1_f(F_v,E_{\pi^n}):= \Image[E(F_v)/\pi^nE(F_v) \to
H^1(F_v,E_{\pi^n})],\\
&H^1_f(F_v,E_{\pi^{*n}}):= \Image[E(F_v)/\pi^{*n}E(F_v) \to
H^1(F_v,E_{\pi^{*n}})].
\end{align*}

Suppose that $M \in \{W,W^*,E_{\pi^n}, E_{\pi^{*n}}\}$ and that
$\fq \in \{\fp,\fp^*\}$ . If $c \in H^1(F,M)$, then we write
$\loc_v(c)$ for the image of $c$ in $H^1(F_v,M)$. We define

$\bullet$ the {\it true Selmer group} $\Sel(F,M)$ by
$$
\Sel(F,M) = \left\{ c \in H^1(F,M) \mid \loc_v(c) \in
H^1_f(F_v,M)\, \text{for all $v$} \right\};
$$

$\bullet$ the {\it relaxed Selmer group} $\Sel_{\rel}(F,M)$ by
$$
\Sel_{\rel}(F,M) = \left\{ c \in H^1(F,M) \mid \loc_v(c) \in
H^1_f(F_v,M)\, \text{for all $v$ not dividing $p$}
\right\};
$$

$\bullet$ the {\it strict Selmer group} $\Sel_{\str}(L,M)$ by
$$
\Sel_{\str}(F,M) = \left\{ c \in \Sel(F,M) \mid \loc_v(c) = 0\,
\text{for all $v$ dividing $p$} \right\};
$$

$\bullet$ the {\it $\fq$-strict Selmer group} $\Sel_{\str(\fq)}(F,M)$
by
$$
\Sel_{\str(\fq)}(F,M)= \left\{ c \in \Sel(F,M) \mid \loc_v(c) = 0\,
\text{for all $v$ dividing $\fq$} \right\};
$$

$\bullet$ the {\it $\fq$-restricted Selmer group} (or simply {\it
restricted Selmer group} for short when $\fq$ is understood)
$\ms_\fq(F,M)$ by
$$
\ms_\fq(F,M) = \left\{ c \in \Sel_{\rel}(F,M) \mid \loc_v(c) = 0\,
\text{for all $v$ dividing $\fq$} \right\}.
$$
(The terminology `restricted Selmer group' is meant to reflect a
choice of a combination of relaxed and strict Selmer conditions at
places above $p$.)

We also define
\begin{align*}
&\ch{\Sel}_{?}(F,T):= \varprojlim_{n} \Sel_{?}(F,E_{\pi^n}),\quad
\ch{\Sel}_{?}(F,T^*):= \varprojlim_{n} \Sel_{?}(F,E_{\pi^{*n}}),\\
&\ch{\ms}_{\fq}(F,T):= \varprojlim_{n} \ms_{\fq}(F,E_{\pi^n}),\quad
\ch{\ms}_{\fq}(F,T^*):= \varprojlim_{n} \ms_{\fq}(F,E_{\pi^{*n}}).
\end{align*}

If $L/K$ is an infinite extension, we define
\begin{align*}
&\Sel_{?}(L,M) = \varinjlim \Sel_{?}(L',M),\quad
\ms_{\fq}(L,M) = \varinjlim \ms_{\fq}(L',M),\\
&\ch{\Sel}_?(L,T)= \varinjlim \ch{\Sel}_?(L',T),\quad 
\ch{\Sel}_?(L,T^*)= \varinjlim \ch{\Sel}_?(L',T^*),
\end{align*}
where the direct limits are taken with respect to restriction over all
subfields $L' \subset L$ finite over $K$.

For any extension $L/K$, we set
$$ 
\Sel_{?}(L,M)^{\land}=X_{?}(L,M), \qquad \ms_\fq(L,M)^{\land}=
X_\fq(L,M).
$$

\begin{theorem} \label{T:coates}
Let $L$ be any field such that $\cF_\infty^* \subseteq L \subseteq
\fF_\infty$. Then there is an isomorphism
\begin{equation} \label{E:coates1}
X_{\fp^*}(L,W^*) \simeq \cX(L)(\psi^{*-1})
\end{equation}
of $\Lambda(L)$-modules.
\end{theorem}

\begin{proof} This is simply the analogue for restricted Selmer groups
of a well-known theorem of Coates concerning true Selmer groups (see
\cite[Theorem 12]{C}). We first observe that, since $\cF_\infty^*
\subseteq L$, we have isomorphisms of $\Lambda(L)$-modules
$$
\cX(L)(\psi^{*-1}) \simeq \Hom(T^*,\cX(L)),\qquad
\cX(L)(\psi^{*-1})^{\land} \simeq \Hom(\cX(L), W^*).
$$
Hence, in order to establish the desired result, it suffices to show
that there is a natural isomorphsim
\begin{equation} \label{E:coates}
\ms_{\fp^*}(L,W^*) \xrightarrow{\sim} \Hom(\cX(L),W^*).
\end{equation}
This may be proved in exactly the same way as \cite[Theorem 12]{C}.
\end{proof}

The following result is a `control theorem' for restricted Selmer
groups.

\begin{proposition} \label{P:control}
(a) Let $I_{\cF_\infty^*}$ denote the kernel of the quotient map
$\Pi_{\cF_\infty^*}: \Lambda(\fF_\infty) \to
\Lambda(\cF_\infty^*)$. Then the kernel of the restriction map
\begin{equation*}
\ms_{\fp^*}(\cF_\infty^*, W^*) \to \ms_{\fp^*}(\fF_\infty,
W^*)[I_{\cF_\infty^*}]
\end{equation*}
is finite. A characteristic power series in
$\Lambda(\cF_{\infty}^{*})$ of the Pontryagin dual of the cokernel of
this map is given by
$$
e_F = (\gamma - \psi^{*-1}(\gamma))^{-1} \prod_{v \mid \fp^*}
(\gamma_v - \psi^{*-1}(\gamma_v)),
$$
where $\gamma$ is a topological generator of
$\Gal(\cF_{\infty}^{*}/F)$, and, for each place $v$ of
$\cF_{\infty}^{*}$ lying above $\fp^*$, $\gamma_v$ denotes a
topological generator of $\Gal(\cF_{\infty,v}^{*}/F_v) \leq
\Gal(\cF_{\infty}^{*}/F)$. 

Hence if $f \in \Lambda(\fF_\infty)$ is a characteristic power series
of $X_{\fp^*}(\cF_\infty^*,W^*)$, then $e_{F}^{-1}
\Pi_{\cF_\infty^*}(f) \in \Lambda(\cF_\infty^*)$ is a characteristic
power series of $X_{\fp^*}(\cF_\infty^*,W^*)$.
\smallskip

(b) Suppose that $L$ is any field such that $F \subseteq L \subseteq
\cF_\infty^*$, and write $I_L$ for the kernel of the quotient map
$\Lambda(\cF_\infty^*) \to \Lambda(L)$. Then the restriction map
\begin{equation*}
\ms_{\fp^*}(L, W^*) \to \ms_{\fp^*}(\cF^*_\infty, W^*)[I_L]
\end{equation*}
is an isomorphism.

Hence the dual of this restriction map is an isomorphism of
$\Lambda(L)$-modules:
\begin{equation*}
X_{\fp^*}(\cF^*_\infty,W^*)/I_LX_{\fp^*}(\cF_\infty,W^*) \xrightarrow{\sim}
X_{\fp^*}(L,W^*).
\end{equation*}
\end{proposition}

\begin{proof} Let $\cN$ denote the maximal extension of $\fF_\infty$
that is unramified away from all places of $\fF_\infty$ lying above
$p$. Consider the following commutative diagram:
\begin{equation*}
\begin{CD}
0 @>>> \ms_{\fp^*}(\cF_\infty^*, W^*) @>>> H^1(\cN/\cF_\infty^*, W^*)
@>{\loc_{\fp^*}}>>\prod_{v \mid \fp^*} H^1(\cN_v/\cF_{\infty,v}^{*}, W^*) \\
@. @V{\alpha}VV  @VVV  @VVV  \\
0 @>>> \ms_{\fp^*}(\fF_\infty, W^*)[I_{\cF_\infty^*}] @>>>
H^1(\cN/\fF_\infty, W^*)[I_{\cF_\infty^*}] @>{\loc_{\fp^*}}>>
\prod_{v \mid \fp^*}H^1(\cN_v/\fF_{\infty,v}, W^*)
\end{CD}
\end{equation*}
in which the vertical arrows are the obvious restriction maps.

Applying the Snake Lemma (together with the inflation-restriction
exact sequence) to this diagram yields the exact sequence
\begin{align} \label{E:snake1}
&0 \to \Ker(\alpha) \to H^1(\fF_\infty/\cF_\infty^*,W^*)
\xrightarrow{g_1} \prod_{v \mid \fp^*}
H^1(\fF_{\infty,v}/\cF_{\infty,v}^{*},W^*) \to \notag \\
&\to \Coker(\alpha) \to H^2(\fF_\infty/\cF_\infty^*,W^*)
\xrightarrow{g_2} \prod_{v \mid \fp^*}
H^2(\fF_{\infty,v}/\cF_{\infty,v}^{*},W^*) \to 0. 
\end{align}

Now,
\begin{align} \label{E:cohom}
&H^1(\fF_\infty/\cF_\infty^*,W^*) \simeq
\Hom(\Gal(\fF_\infty/\cF_\infty^*),W^*), \notag \\
&\prod_{v \mid \fp^*} H^1(\fF_{\infty,v}/\cF_{\infty,v}^{*},W^*) \simeq
\prod_{v \mid \fp^*} \Hom(\Gal(\fF_{\infty,v}/\cF_{\infty,v}^{*}),W^*),
\end{align}
and, as $\Gal(\fF_\infty/\cF_\infty^*) \simeq \Delta \times \bZ_p$
with $p \nmid \Delta$, we have
\begin{align*}
&H^2(\fF_\infty/\cF_\infty^*,W^*) \simeq
H^0(\fF_\infty/\cF_\infty^*,W^*) \simeq W^*, \\
&\prod_{v \mid \fp^*} H^2(\fF_{\infty,v}/\cF_{\infty,v}^{*},W^*)
\simeq
\prod_{v \mid \fp^*} H^0(\fF_{\infty,v}/\cF_{\infty,v}^{*},W^*)
\simeq
\prod_{v \mid \fp^*} W^*.
\end{align*}

We now deduce that $g_1$ is non-zero, and therefore has finite kernel
(since $H^1(\fF_\infty/\cF_{\infty}^{*},W^*)$ is divisible), and that
$g_2$ is injective. It follows from \eqref{E:snake1} that
$\Ker(\alpha)$ is finite, and that there is an exact sequence
\begin{equation} \label{E:snake2}
0 \to \Ker(\alpha) \to H^1(\fF_\infty/\cF_\infty^*,W^*)
\xrightarrow{g_1} \prod_{v \mid \fp^*}
H^1(\fF_{\infty,v}/\cF_{\infty,v}^{*},W^*)
\to \Coker(\alpha) \to 0.
\end{equation}

It follows from \eqref{E:cohom} that
\begin{align*}
&\cha_{\Lambda(\cF_{\infty}^{*})} \left(
H^1(\fF_\infty/\cF_\infty^*,W^*) \right)^{\land} 
= \gamma - \psi^{*-1}(\gamma); \\
&\cha_{\Lambda(\cF_\infty^*)} \left( \prod_{v \mid \fp^*}
  H^1(\fF_{\infty,v}/\cF_{\infty,v}^{*}, W^*) \right)^{\land} =
\prod_{v \mid \fp^*} (\gamma_v - \psi^{*-1}(\gamma_v)).
\end{align*}
Hence we deduce from \eqref{E:snake2} that 
$$
\cha_{\Lambda(\cF_{\infty}^{*})}(\Coker(\alpha))^{\land} = e_F = 
(\gamma - \psi^{*-1}(\gamma))^{-1} \prod_{v \mid \fp^*}
(\gamma_v - \psi^{*-1}(\gamma_v)),
$$
as asserted.
\smallskip

(b) In this case we consider the commutative diagram
\begin{equation*}
\begin{CD}
0 @>>> \ms_{\fp^*}(L, W^*) @>>> H^1(\cN/L, W^*)
@>{\loc_{\fp^*}}>>\prod_{v \mid \fp^*} H^1(\cN_v/L_{v}, W^*) \\
@. @V{\beta_1}VV  @V{\beta_2}VV  @V{\beta_3}VV  \\
0 @>>> \ms_{\fp^*}(\cF_\infty^*, W^*)[I_L] @>>>
H^1(\cN/\cF^*_\infty, W^*) @>{\loc_{\fp^*}}>>
\prod_{v \mid \fp^*}H^1(\cN/\cF^{*}_{\infty,v}, W^*)
\end{CD}
\end{equation*}

We have that
\begin{align*}
&\Ker(\beta_2) = H^1(\cF_\infty^*/L,W^*) = 0, \\
&\Ker(\beta_3) = \prod_{v \mid \fp^*} H^1(\cF_{\infty,v}^{*}/L_v,W^*)
=0, \\
&\Coker(\beta_2) = H^2(\cF_\infty^*/L,W^*) = 0,
\end{align*}
(see \cite[p. 40]{PR7}, for example), and so the Snake Lemma implies
that $\beta_1$ is an isomorphism, as claimed.
\end{proof}

\begin{corollary} \label{C:descent}
For any field $L$ with $F \subseteq L \subseteq \cF^*_\infty$, we have
an isomorphism
\begin{equation} \label{E:descent}
X_{\fp^*}(L,T^*) \simeq
\cX(\cF_\infty)(\psi^{*-1})/I_L(\cX(\cF_\infty)(\psi^{*-1})
\end{equation}
of $\Lambda(L)$-modules.
\end{corollary}

\begin{proof} This follows directly from Proposition \ref{P:control}
and Theorem \ref{T:coates}.
\end{proof}

\begin{remark} \label{R:char} 
If we take $F=K$ in Proposition \ref{P:control}, then it is easy to
check that $e_K \in \Lambda(\cK_{\infty}^{*})^{\times}$.  We therefore
see from Proposition \ref{P:control}(a) and Corollary \ref{C:descent}
that the element $H_K \in \Lambda(K^*_\infty)$ fixed in Section
\ref{S:twists} is a characteristic power series of
$X_{\fp^*}(K^*_\infty,W^*)$. \qed
\end{remark}

\begin{definition} \label{D:sha2}
For any finite extension $F/K$ and any prime $\fq$ of $K$ we define
$$
\sha(F)_{\rel(\fq)}:= \Ker \left[ H^1(F,E) \to \prod_{v \nmid \fq}
  H^1(F_v,E) \right],
$$
and we set
$$
E_{1,\fq}(F):= \Ker \left[ E(F) \otimes_{O_K} O_{K,\fq} \to \prod_{v \mid
  \fq} E(F_v) \right].
$$
\qed
\end{definition}

\begin{lemma} \label{L:free}
Let $F/K$ be any finite extension, and let $\fq \in \{\fp, \fp^*\}$.
Then $\ch{\ms}_{\fq}(F,T_{\fq})$ is a free $O_{K,\fq}$-module.
\end{lemma}

\begin{proof} 
It follows from the definitions that
$\ch{\ms}_{\fq}(F,T_{\fq})_{\tors} \subseteq
\ch{\Sel}(F,T_{\fq})$. The desired result now follows from the fact
that the restriction of the localisation map
$$
\ch{\Sel}(F,T_{\fq}) \to \prod_{v \mid \fq} E(F_v) \otimes_{O_K}
O_{K,\fq}
$$
to $\ch{\Sel}(F,T_{\fq})_{\tors}$ is injective.
\end{proof}

\section{The $p$-adic height pairing on restricted Selmer groups}
\label{S:height} 

In this section we shall explain how the methods described by
Perrin-Riou in \cite{PR1} and \cite{ PR7} may be used to construct a
$p$-adic height pairing
$$
[\,, \,]_{F,\fp^*}: \ms_\fp(F,T) \times \ms_{\fp^*}(F,T^*) \to O_{K,\fp^*}.
$$

We begin by describing the $\fp$-adic Leopoldt hypotheses with which we
shall work.

\begin{definition} \label{D:leo}
Let $M/K$ be any finite extension, and consider the diagonal injection
$$
i_M: O_M^\times \to \prod_{v \mid \fp} O_{M,v}^{\times}.
$$
Let $\ov{i_M(O_M^\times)}$ denote the $\fp$-adic closure of
$i_M(O_M^\times)$ in $\prod_{v \mid \fp} O_{M,v}^{\times}$, and set
$$
\delta(M):= \rk_{\bZ}(O_M^\times) - \rk_{\bZ_p}(\ov{i_M(O_M^\times)}).
$$
The {\it weak $\fp$-adic Leopoldt hypothesis for $F$} asserts that
the numbers $\delta(L')$ are bounded as $L'$ runs through all finite
extensions of $F$ contained in $\cF_\infty^*$. The {\it strong
$\fp$-adic Leopoldt hypothesis for $F$} asserts that the numbers
$\delta(L')$ are all equal to zero. 

We remark that the strong Leopoldt hypothesis is known to hold for all
abelian extensions of $K$ (see \cite{B}). \qed
\end{definition}

Recall that $\cB(\cF^*_\infty)$ denotes the maximal abelian pro-$p$
extension of $\cF_\infty^*$ which is unramified away from $\fp$ and
totally split at all places above $\fp^*$, and that $\cY(\cF^*_\infty)
= \Gal(\cB(\cF^*_\infty)/\cF_\infty^*)$. The main ingredient in the
construction of $[\,,\,]_{F,\fp^*}$ is the following result.

\begin{theorem} \label{T:caniso}
If the weak $\fp$-adic Leopoldt hypothesis holds for $F$ then there
is a natural isomorphism
$$
\Psi_F: \ch{\ms}_{\fp}(F,T) \xrightarrow{\sim}
\Hom(T^*,\cY(\cF_\infty^*))^{\Gal(\cF_\infty^*/F)}.
$$
\end{theorem}

The proof of this theorem is very similar to that of \cite[Th\'eor\`eme
3.2]{PR1}. We shall therefore just describe the main outlines of
the proof, and we refer the reader to \cite{PR1} for some of the
details which we omit.

In order to describe the proof of Theorem \ref{T:caniso}, we require a
number of intermediary results.

\begin{lemma} \label{L:kummer1}
There is an isomorphism of $\Gal(\cF^*_n/F)$-modules
\begin{equation} \label{E:mystiso}
H^1(\cF^*_n,E_{\pi^{n}}) \xrightarrow{\sim}
\Hom(E_{\pi^{*n}},\cF_n^{*\times}/\cF_{n}^{*\times p^n});\quad f \mapsto
\tilde{f}.
\end{equation}
For each place $v$ of $\cF^*_n$, there is also a corresponding local
isomorphism
$$
H^1(\cF^{*}_{n,v},E_{\pi^{n}}) \xrightarrow{\sim}
\Hom(E_{\pi^{*n}},\cF_{n,v}^{*\times}/\cF_{n,v}^{*\times p^n}).
$$
\end{lemma}

\begin{proof} See \cite[Lemme 3.8]{PR1}. The isomorphism
\eqref{E:mystiso} is defined as follows. Let $f \in
H^1(\cF_n^*,E_{\pi^{n}})$, and write
$$
w_n: E_{\pi^n} \times E_{\pi^{*n}} \to \mu_{p^n}
$$ 
for the Weil pairing. We identify $\cF_n^{*\times}/\cF_{n}^{*\times
p^n}$ with $H^1(\cF_n^*, \mu_{p^n})$ via Kummer theory. If $u \in
E_{\pi^{*n}}$, then $\tilde{f}(u) \in H^1(\cF_n^*, \mu_{p^n})$ is
defined to be the element represented by the cocycle
$$
\sigma \mapsto w_n(f(\sigma),u)
$$
for all $\sigma \in \Gal(\ov{F}/\cF_n^*)$.
\end{proof}

\begin{lemma} \label{L:kummer2}
For each place $v$ of $\cF^*_n$ with $v \nmid \fp^*$, there is
an isomorphism
$$
E(\cF^{*}_{n,v})/\pi^{n}E(\cF^{*}_{n,v}) \xrightarrow{\sim}
\Hom(E_{\pi{*^n}}, O_{\cF^{*}_{n,v}}^{\times}/O_{\cF^{*}_{n,v}}^{\times p^n}).
$$
\end{lemma}

\begin{proof} See \cite[Lemme 3.11]{PR1}.
\end{proof}

\begin{corollary} \label{C:kummer3}
Suppose that $h \in H^1(\cF^*_n, E_{\pi^{n}})$. Then $h \in \ms_{\fp}(\cF^*_n,
E_{\pi^{n}})$ if and only if, for each $u \in E_{\pi^{n}}$, the
following local conditions are satisfied:

(a) $\tilde{h}(u) \in \cF_{n,v}^{*\times p^n}$ for all $v \mid \fp$;

(b) $p^n \mid v_{\cF^*_n}(\tilde{h}(u))$ for all $v \nmid \fp^*$.

(Note that we impose no local conditions at places lying above $\fp^*$.)
\end{corollary}

\begin{proof} This follows directly from Lemmas \ref{L:kummer1} and
\ref{L:kummer2}.
\end{proof}

In what follows, we set $G_n:= \Gal(\cF^*_n/F)$, and we write $J_n$ for
the group of finite ideles of $\cF^*_n$. We let $V_n$ denote the subgroup
of $J_n$ consisting of those elements whose components are equal to
$1$ at all places dividing $\fp$ and are units at all places not
dividing $\fp^*$. We set
$$
C_n:= J_n/V_n \cF_{n}^{*\times},\qquad \Omega_n:= \prod_{v \mid \fp}
\mu_{p^n}(\cF^{*}_{n,v}),
$$
and we note that the order of $\Omega_n$ is bounded as $n$ varies.

\begin{proposition}
There is an exact sequence
$$
\Hom(E_{\pi^{*n}},\Omega_n)^{G_n} \to \Hom(E_{\pi^{*n}},C_n)^{G_n}
\xrightarrow{\eta_n} \ms_{\fp}(F,E_{\pi^{n}}) \to 0.
$$
\end{proposition}

\begin{proof} The proof of this Proposition is identical,
\textit{mutatis mutandis}, to that of \cite[Proposition 3.13]{PR1}.
\end{proof}

Now let $\eta'_n$ be the map obtained from $\eta_n$ via passage to the
quotient by the kernel of $\eta_n$, and write $C_n(p)$ for the
$p$-primary part of $C_n$. Then it may be shown exactly as on
\cite[pp. 387--389]{PR1} that passing to inverse limits over the maps
$\eta'^{-1}_{n}$ yields an isomorphism
$$
\Xi_F: \varprojlim \ch{\ms}_{\fp}(F,E_{\pi^{n}}) = \ms_{\fp}(F,T)
\xrightarrow{\sim} \Hom(T^*, \varprojlim C_n(p))^{\Gal(\cF^*_\infty/F)}.
$$
(Here the inverse limit $\varprojlim C_n(p)$ is taken with respect to
the norm maps $\cF_{n}^{*\times} \to \cF_{n-1}^{*\times}$.)

The proof of Theorem \ref{T:caniso} is completed by the following
result.

\begin{proposition} If the weak $\fp$-adic Leopoldt hypothesis holds
for $F$, then there is an isomorphism
$$
\Hom(T^*, \varprojlim C_n(p))^{\Gal(\cF^*_\infty/F)} \simeq
\Hom(T^*,\cY(\cF_\infty^*))^{\Gal(\cF_\infty^*/F)}.
$$
\end{proposition}

\begin{proof} This may be shown in the same way as \cite[Lemme 3.18]{PR1}.
\end{proof}

We now explain how the isomorphism $\Psi_F$ may be used to construct a
$p$-adic height pairing
$$
[\,,\,]_{F,\fp^*}: \ch{\ms}_\fp(F,T) \times \ch{\ms}_{\fp^*}(F,T^*)
\to O_{K,\fp^*}.
$$

We first recall (see Proposition \ref{P:control}(b)) that the restriction
map
\begin{equation} \label{E:res}
\ms_{\fp^*}(F,W^*) \to \ms_{\fp^*}(\cF^*_\infty,W^*)
\end{equation}
is injective, and that there is a natural isomorphism (see Theorem
\ref{T:coates})
\begin{equation} \label{E:coates2}
\ms_{\fp^*}(\cF^*_\infty, W^*) \xrightarrow{\sim}
\Hom(\cX(\cF^*_\infty), W^*).
\end{equation}
It follows from the local conditions defining the restricted Selmer group
$\ms_{\fp^*}(F,W^*)$ that \eqref{E:res} and \eqref{E:coates2} induce
an injection
\begin{equation} \label{E:msgal1}
\ms_{\fp^*}(F,W^*) \rightarrow \Hom(\cY(\cF^*_\infty), W^*),
\end{equation}
and taking Pontryagin duals yields a surjection
\begin{equation} \label{E:msgal2}
\Hom(T^*, \cY(\cF^*_\infty)) \rightarrow X_{\fp^*}(F,W^*).
\end{equation}
Composing this with the natural surjection
$$
X_{\fp^*}(F,W^*) \rightarrow [\ms_{\fp^*}(F,W^*)_{\Div}]^\land
$$
and taking $\Gal(\cF^*_\infty/F)$-invariants yields a homomorphism
$$
\beta_F: \Hom(T^*, \cY(\cF^*_\infty))^{\Gal(\cF^*_\infty/F)}
\rightarrow [\ms_{\fp^*}(F,W^*)_{\Div}]^{\land}.
$$
Next, we observe that we have a canonical isomorphism
\begin{align*}
[\ms_{\fp^*}(F,W^*)_{\Div}]^{\land} &\simeq
\Hom_{O_{K,\fp^*}}(T_{\fp^*}(\ms_{\fp^*}(F,W^*)_{\Div}),O_{K,\fp^*}) \\
&= \Hom_{O_{K,\fp^*}}(T_{\fp^*}(\ms_{\fp^*}(F,W^*)) , O_{K,\fp^*}),
\end{align*}
where the last equality holds because
$$
T_{\fp^*}(\ms_{\fp^*}(F,W^*)_{\Div} = T_{\fp^*}(\ms_{\fp^*}(F,W^*)).
$$
Also, for each $n \geq 1$, we have a surjective map
$$
\ms_{\fp^*}(F,E_{\pi^{*n}}) \to \ms_{\fp^*}(F,W^*)_{\pi^{*n}}
$$
with finite kernel. Via passage to inverse limits, these yield a
map
$$
\ch{\ms}_{\fp^*}(F,T^*) \to T_{\fp^*}(\ms_{\fp^*}(F,W^*))
$$ 
which is an isomorphism because $\ch{\ms}_{\fp^*}(F,T^*)$ is
$O_{K,\fp^*}$-free (see Lemma \ref{L:free}).

It follows from the above discussion that we may view $\beta_F$ as a
homomorphism
$$
\beta_F: \Hom(T^*, \cY(\cF^*_\infty))^{\Gal(\cF^*_\infty/F)} \to
\Hom_{O_{K,\fp^*}}(\ch{\ms}_{\fp^*}(F,T^*) , O_{K,\fp^*}).
$$
We thus obtain a map
$$
\beta_F \circ \Psi_F: \ch{\ms}_{\fp}(F,T) \to
\Hom_{O_{K,\fp^*}}(\ch{\ms}_{\fp^*}(F,T^*) , O_{K,\fp^*}),
$$
and this yields the desired pairing
$$
[\,,\,]_{F,\fp^*}: \ch{\ms}_\fp(F,T) \times \ch{\ms}_{\fp^*}(F,T^*)
\to O_{K,\fp^*}.
$$ 
It is natural to conjecture that this pairing is always
non-degenerate (see Remark \ref{R:height}).

If $x_1,\ldots,x_m$ is an $O_{K,\fp}$-basis of $\ch{\ms}_{\fp}(F,T)$
(resp. if $y_1,\ldots , y_m$ is an $O_{K,\fp^*}$-basis
of $\ch{\ms}_{\fp^*}(F,T^*)$), then we define the regulator
$\cR_{F,\fp^*}$ associated to $[\,,\,]_{F,\fp^*}$ by
\begin{equation} \label{E:reg}
\cR_{F,\fp^*}:= \det([x_i,y_j]_{F,\fp^*}).
\end{equation}

\section{The leading term} \label{S:leading}

We retain the notation of the previous section. Write $\Gamma_F:=
\Gal(F_\infty^*/F)$, fix a topological generator $\gamma_F$ of
$\Gamma_F$, and identify $\Lambda(F_\infty^*)$ with the power series
ring $\bZ_p[[t]]$ via the map $\gamma_F \mapsto t+1$. Let $H_F \in
\Lambda(F^*_\infty)$ be a characteristic power series of
$X_{\fp^*}(F_\infty^*,W^*)$. In this section we shall calculate the
leading coefficient of $H_F$, assuming that the strong Leopoldt
hypothesis holds for $F$ and that $[\,,\,]_{F,\fp^*}$ is
non-degenerate.

\begin{proposition} \label{P:nofinite}
Suppose that $F$ satisfies the strong $\fp$-adic Leopoldt
hypothesis. Then the $\Lambda(F_\infty^*)$-module
$X_{\fp^*}(F_\infty^*,W^*)$ has no finite, non-trivial submodules.
\end{proposition}

\begin{proof} 
It is straightforward to show that a slight modification of the
arguments given in \cite[\S4]{Gr78} establishes the fact that if $F$
satisfies the strong $\fp$-adic Leopoldt hypothesis, then the
$\Lambda(F_\infty^*)$-module $X(F_\infty^*)$ has no finite,
non-trivial submodules. For brevity, we omit the details. The desired
result now follows from Proposition \ref{P:control} and Theorem
\ref{T:coates}.
\end{proof}

\begin{theorem} \label{T:leading}
Let $H_F \in \Lambda(F^*_\infty)$ be a characteristic power series of
$X_{\fp^*}(F_\infty^*,W^*)$. Assume that the strong $\fp$-adic
Leopoldt hypothesis holds for $F$, and that $[\,,\,]_{F,\fp^*}$ is
non-degenerate. Set $m:=
\rk_{O_{K,\fp^*}}(\ch{\ms}_{\fp^*}(F,T^*))$. Then $\ord_{t=0} H_F =m$,
and
\begin{equation} \label{E:leading}
\frac{H_F}{t^m} \Biggr |_{t=0} \sim |\ms_{\fp^*}(F,W^*)_{/\Div}| \cdot
\cR_{F,\fp^*}.
\end{equation}
\end{theorem}

\begin{proof}
We begin by noting that there is a surjective homomorphism
$$
X_{\fp^*}(F_\infty^*, W^*) \rightarrow
[\ms_{\fp^*}(F,W^*)_{\Div}]^{\land}.
$$
This implies that $H_F$ is divisible by $t^m$. If we write $Z_\infty$
for the kernel of this map, then the Snake Lemma yields the following
exact sequence:
\begin{align*}
0 &\to (Z_\infty)^{\Gamma_F} \to X_{\fp^*}(F_\infty^*,W^*)^{\Gamma_F}
\xrightarrow{\xi_F} [\ms_{\fp^*}(F,W^*)_{\Div}]^{\land} \to \\
&\to (Z_\infty)_{\Gamma_F} \to X_{\fp^*}(F_\infty^*, W^*)_{\Gamma_F}
\to [\ms_{\fp^*}(F,W^*)_{\Div}]^{\land} \to 0.
\end{align*}

The kernel of the last map
$$
X_{\fp^*}(F_\infty^*, W^*)_{\Gamma_F} \to
[\ms_{\fp^*}(F,W^*)_{\Div}]^{\land}
$$
is dual to the cokernel of the map
$$
\ms_{\fp^*}(F,W^*)_{\Div} \to \ms_{\fp^*}(F^*_\infty, W^*)^{\Gamma_F}.
$$ 
Since $\ms_{\fp^*}(F,W^*) \simeq \ms_{\fp^*}(F^*_\infty,
W^*)^{\Gamma_F}$ (via Proposition \ref{P:control}(b)), it follows that
this cokernel is isomorphic to $\ms_{\fp^*}(F,W^*)_{/\Div}$, which is
finite.

We therefore deduce that the multiplicity of $t$ in $H_F$ is equal to
$m$ if and only if $(Z_\infty)_{\Gamma_F}$ is finite, which in turn is
the case if and only if the cokernel of $\xi_F$ is finite. Recall (see
Theorem \ref{T:coates})
$$
X_{\fp^*}(F_\infty^*,W^*)^{\Gamma_F} \simeq \Hom(T^*,
\cX(\cF_\infty^*))^{\Gal(\cF_\infty^*/F)},
$$
and that the homomorphism $\xi_F$ may be written as the following
composition of maps
$$
\Hom(T^*,\cX(F^*_\infty))^{\Gal(\cF^*_\infty/F)} \to
\Hom(T^*,\cY(F^*_\infty))^{\Gal(\cF^*_\infty/F)} \to
\ms_{\fp^*}(F,W^*)^{\land} \to 
[\ms_{\fp^*}(F,W^*)_{/\Div}]^{\land}
$$
(see \eqref{E:msgal1}, \eqref{E:msgal2}). Hence the cokernel of
$\xi_F$ is finite if and only if the $p$-adic height pairing
$[\,,\,]_{F,\fp^*}$ is non-degenerate.

We now see that if $[\,,\,]_{F,\fp^*}$ is non-degenerate, then
$(Z_{\infty})_{\Gamma_F}$ is finite. This implies that
$(Z_\infty)^{\Gamma_F}$ is also finite, whence it follows via
Proposition \ref{P:nofinite} that $(Z_\infty)^{\Gamma_F}=0$. Hence we
have
$$
\frac{H_F}{t^m} \Biggr |_{t=0} \sim |(Z_{\infty})_{\Gamma_F}| \sim
|\ms_{\fp^*}(F,W^*)_{/\Div}| \cdot |\Coker(\xi_F)|.
$$
Now
\begin{align*}
|\Coker(\xi_F)| &= [(\ms_{\fp^*}(F,W^*)_{\Div})^{\land}:
 \xi_F(X_{\fp^*}(F_\infty^*, W^*)^{\Gamma_F})] \\
&= [T_{\fp^*}(\ms_{\fp^*}(F,W^*)): \Psi_F(\ch{\ms}_{\fp}(F,T))] \\
&= \cR_{F,\fp^*} \cdot \left[ \Ker(\ch{\ms}_{\fp^*}(F,T^*) \to
 T_{\fp^*}(\ms_{\fp^*}(F,W^*))) \right] \\
&= \cR_{F,\fp^*}.
\end{align*}

Hence
$$
\frac{H_F}{t^m} \Biggr |_{t=0} \sim |\ms_{\fp^*}(F,W^*)_{/\Div}| \cdot
\cR_{F,\fp^*},
$$
as claimed.
\end{proof}

\section{Restricted Selmer groups over $K$} \label{S:rsk}

In this section we shall analyse various properties of restricted
Selmer groups over $K$. The main tool for doing this is the
Poitou-Tate exact sequence (see e.g. \cite[Theorem 1.5]{CS} or
\cite[Proposition 4.1.1]{PR6}).

We write $S_F$ for the set of places of $F$ lying above $p$, and
$G_{F,S_F}$ for the Galois group over $F$ of the maximal abelian
extension of $F$ that is unramified away from all places in $S_F$.

\begin{proposition} \label{P:costr}
There are isomorphisms
$$
\ch{\Sel}_{\str}(F,T^*)
\simeq H^2(G_{F,S_F},W)^{\land}, \qquad
\ch{\Sel}_{\str}(F,T) \simeq H^2(G_{F,S_F},W^*)^{\land}.
$$
\end{proposition}

\begin{proof} The middle of the Poitou-Tate exact sequence yields
$$
0 \to \Sel_{\str}(F,E_{\pi^{*n}})^{\land} \to H^2(G_{F,S_F}, E_{\pi^n}) \to
\bigoplus_{v \in S_F} H^2(F_v, E_{\pi^n}).
$$
Dualising, and using the fact that, via Tate local duality, we have
$H^2(F_v,E_{\pi^n})^{\land} \simeq H^0(F_v, E_{\pi^{*n}})$ for each
place $v$ of $F$ gives
$$
\bigoplus_{v \in S_F} H^0(F_v, E_{\pi^{*n}}) \to H^2(G_{F,S_F},
E_{\pi^n})^{\land} \to \Sel_{\str}(F,E_{\pi^{*n}}) \to 0.
$$
By passing to limits we obtain
$$
\bigoplus_{v \in S_F} H^0(F_v, T^*) \to H^2(G_{F,S_F},W)^{\land} \to
\ch{\Sel}_{\str}(F,T^*) \to 0,
$$
and this establishes the first isomorphism, since the first term of
this last sequence is equal to zero.

The second isomorphism may be proved in a similar manner.
\end{proof}

Recall that $r= \rk_{O_K}(E(K))$.

\begin{proposition} \label{P:rkstr1}
Suppose that $r \geq 1$. Then
\begin{align*}
\rk_{O_{K,\fp^*}}(\ch{\Sel}_{\str}(K,T^*)) &=
\rk_{O_{K,\fp^*}}(\ch{\Sel}_{\str(\fp^*)}(K,T^*)) \\
&= \rk_{O_{K,\fp^*}}(\ch{\Sel}(K,T^*))
- 1.
\end{align*}
\end{proposition}

\begin{proof}
Since $r \geq 1$, the image of the localisation map
$$
\Sel(K,T^*) \to E(K_{\fp^*}) \otimes O_{K,\fp^*}
$$
is infinite. The result now follows from the fact that
$$
\rk_{O_{K,\fp^*}}[E(K_{\fp^*}) \otimes O_{K,\fp^*}]=
\rk_{O_{K,\fp^*}}\left[ \prod_{v \mid p} E(K_v) \otimes
O_{K,\fp^*}\right] =1.
$$
\end{proof}

\begin{lemma} \label{L:Euler}
(a) The cohomology group $H^1_f(K_{\fp^*},T)$ is finite, and
$$
| H^1_f(K_{\fp^*},T) | \sim |\tilde{E}_{\fp^*}(k_{\fp^*}) |
\sim 1 - \psi(\fp^*)
$$
in $\bZ_{p}$.

(b) We have
$$
H^1_f(K_{\fp^*},T) = H^1(K_{\fp^*},T)_{\tors},
$$
and $H^1(K_{\fp^*},T)/H^1_f(K_{\fp^*},T)$ is $O_{K,\fp^*}$-free of
rank one.
\end{lemma}

\begin{proof} 
Part (a) follows directly from \cite[Lemma 1]{C}.

To prove part (b), we observe that, via Tate local
duality, the dual of $H^1(K_{\fp^*},T)/H^1_f(K_{\fp^*},T)$ is equal to
$E(K_{\fp^*}) \otimes D_{\fp^*}$, and this last group is divisible of
$O_{K,\fp^*}$-corank one.
\end{proof}

\begin{proposition} \label{P:rkstr2}
(a) Suppose that $r \geq 1$. Then
$$
\rk_{O_{K,\fp^*}}(\ch{\Sel}_{\rel}(K,T^*)) =
\rk_{O_{K,\fp^*}}(\ch{\Sel}(K,T^*)),
$$
and
$$
[\ch{\Sel}_{\rel}(K,T^*): \ch{\Sel}(K,T^*)] \sim
|\tilde{E}_{\fp^*}(k_{\fp^*})|. 
$$

(b) Suppose that $r=0$. Then
$$
\rk_{O_{K,\fp^*}}(\ch{\Sel}_{\rel}(K,T^*)) = 1.
$$
\end{proposition}

\begin{proof} The Poitou-Tate exact sequence yields
\begin{equation} \label{E:comp1}
0 \to \ch{\Sel}(K,T^*) \to \ch{\Sel}_{\rel}(K,T^*) \xrightarrow{\alpha}
\bigoplus_{v \mid p} \frac{H^1(K_v,T^*)}{H^1_f(K_v,T^*)} \to
\Sel(K,W)^{\land}.
\end{equation}
The cokernel of $\alpha$ is 
the Pontryagin dual of the image
of the localisation map
$$
\Sel(K,W) \to \bigoplus_{v \mid p} H^1_f(K_v,W),
$$
and so has $O_{K,\fp^*}$-rank one if $r \geq 1$ and rank zero if
$r=0$. As
$$
\rk_{O_{K,\fp^*}}[ \oplus_{v \mid
p}(H^1(K_v,T^*)/H^1_f(K_v,T^*))] = 1,
$$
we therefore deduce that $\rk_{O_{K,\fp^*}}(\ch{\Sel}_{\rel}(K,T^*))$
is equal to $\rk_{O_{K,\fp^*}}(\ch{\Sel}(K,T^*))$ if $r \geq 1$, and
is equal to one if $r =0$. In particular, we have that
$\ch{\Sel}_{\rel}(K,T^*)/\ch{\Sel}(K,T^*)$ is finite if $r \geq 1$.

Now suppose that $r \geq 1$. As $H^1(K_{\fp},T^*)/H^1_f(K_\fp,T^*)$ is
$O_{K,\fp^*}$-free of rank one (Lemma \ref{L:Euler}(b)) and
$\ch{\Sel}_{\rel}(K,T^*)/\ch{\Sel}(K,T^*)$ is finite, \eqref{E:comp1}
implies that there is an exact sequence   
$$
0 \to \frac{\ch{\Sel}_{\rel}(K,T^*)}{\ch{\Sel}(K,T^*)} \to
\frac{H^1(K_{\fp^*},T^*)}{H^1_f(K_{\fp^*},T^*)} \xrightarrow{\alpha'} 
\Sel(K,W)^{\land}.
$$
Since $E(K_{\fp^*}) \otimes D_{\fp}=0$, it follows that $\alpha'$ is
the zero map. The dual of $H^1(K_{\fp^*},T^*)/H^1_f(K_{\fp^*},T^*)$ is
isomorphic to $H^1_f(K_{\fp^*},T)$, and Lemma \ref{L:Euler}(a) implies
that
$$
|H^1_f(K_{\fp^*},T)| \sim | \tilde{E}_{\fp^*}(k_{\fp^*})|.
$$
Hence $[\ch{\Sel}_{\rel}(K,T^*): \ch{\Sel}(K,T^*)] \sim
|\tilde{E}_{\fp^*}(k_{\fp^*})|$, as claimed.
\end{proof}

\begin{proposition} \label{P:rkms}
 Suppose that $r \geq 1$. Then
$$
\ch{\ms}_{\fp^*}(K,T^*) = \ch{\Sel}_{\str(\fp^*)}(K,T^*).
$$
In particular, we have
$$
\rk_{O_{K,\fp^*}}(\ch{\ms}_{\fp^*}(K,T^*)) =
\rk_{O_{K,\fp^*}}(\ch{\Sel}(K,T^*)) - 1. 
$$
\end{proposition}

\begin{proof} From Proposition \ref{P:rkstr2}(a), we have
$$
\rk_{O_{K,\fp^*}}(\ch{\Sel}_{\rel}(K,T^*)) =
\rk_{O_{K,\fp^*}}(\ch{\Sel}(K,T^*)).
$$
This implies that
\begin{align} \label{E:rkms}
\rk_{O_{K,\fp^*}}(\ms_{\fp^*}(K,T^*)) &=
\rk_{O_{K,\fp^*}}(\ch{\Sel}_{\str(\fp^*)}(K,T^*)) \notag \\
&= \rk_{O_{K,\fp^*}}(\ch{\Sel}(K,T^*)) - 1. 
\end{align}

It follows from the definitions of $\ch{\ms}_{\fp^*}(K,T^*)$ and
$\ch{\Sel}_{\str(\fp^*)}(K,T^*)$ that we have the following exact
sequence
$$
0 \to \ch{\Sel}_{\str(\fp^*)}(K,T^*) \to \ch{\ms}_{\fp^*}(K,T^*)
\xrightarrow{\beta} \frac{H^1(K_{\fp^*},T^*)}{H^1_f(K_{\fp},T^*)} \to
\Coker(\beta) \to 0,
$$
where $\beta$ is induced by the obvious localisation map. From
\eqref{E:rkms}, we see that
$\ch{\ms}_{\fp^*}(K,T^*)/\ch{\Sel}_{\str(\fp^*)}(K,T^*)$ is
finite. Hence, as $H^1(K_{\fp},T^*)/H^1_f(K_{\fp},T^*)$ is
$O_{K,\fp^*}$-free of rank one (see Lemma \ref{L:Euler}(b)), it
follows that $\beta$ is the zero map. This implies that
$$
\ch{\ms}_{\fp^*}(K,T^*) = \ch{\Sel}_{\str(\fp^*)}(K,T^*)
$$
as claimed.

The final assertion of the Proposition is a direct consequence of
Proposition \ref{P:rkstr1}.
\end{proof}

\begin{remark} \label{R:height} 
Suppose that $r \geq 1$. Then it follows from Proposition
\ref{P:rkms}, together with the definition of $[\,,\,]_{K,\fp^*}$ that
the pairing $[\,,\,]_{K,\fp^*}$ is simply the restriction of
Perrin-Riou's algebraic $p$-adic height pairing $\{\,,\,\}_{K,\fp^*}$
to $\ch{\Sel}_{\str(\fp^*)}(K,T^*) \times
\ch{\Sel}_{\str(\fp)}(K,T)$. Hence, if $r \geq 1$ and
$\{\,,\,\}_{K,\fp^*}$ is non-degenerate, then so is
$[\,,\,]_{K,\fp^*}$. We conjecture that the pairing
$[\,,\,]_{K,\fp^*}$ is also non-degenerate when $r=0$. \qed
\end{remark}

\begin{proposition} \label{P:rkms2}
Suppose that $r=0$. Then
$$
\rk_{O_{K,\fp^*}}(\ch{\ms}_{\fp^*}(K,T^*)) = 1.
$$
\end{proposition}

\begin{proof} We have an injection
$$
0 \to \ch{\ms}_{\fp^*}(K,T^*) \to \ch{\Sel}_{\rel}(K,T^*),
$$
and we know that $\rk_{O_{K,\fp^*}}(\ch{\Sel}_{\rel}(K,T^*))=1$
(Proposition \ref{P:rkstr2}(b)). Hence
$\rk_{O_{K,\fp^*}}(\ch{\ms}_{\fp^*}(K,T^*))$ is either zero or one.

Suppose that $\rk_{O_{K,\fp^*}}(\ch{\ms}_{\fp^*}(K,T^*)) = 0.$ Then
the proof of Theorem \ref{T:leading} shows that the characteristic
power series $H_K \in \Lambda(K_\infty^*)$ of $X_{\fp^*}(K,W^*)$ does
not vanish at $t=0$. This implies that $\ord_{s=1}L_{\fp}^{*}(s) =0$
(see \eqref{E:orders}). On the other hand, it follows from the
functional equation satisfied by the two-variable $p$-adic
$L$-function $\cL_{\fp}$ (see \cite[Chapter II, \S6]{dS}) that the
orders of the zeros at $s=1$ of $L_{\fp}(s)$ and $L_{\fp^*}(s)$ have
opposite parity. Since $r=0$, the order of $\sha(K)$ is known to be
finite (see \cite{Ru}), and so
$$
\ord_{s=1} L_{\fp}(s) = \rk_{O_{K,\fp^*}}(\Sel(K,T^*)) = 0.
$$
This implies that $\ord_{s=1}L_{\fp}^{*}(s) \geq 1$, which is a
contradiction.

It therefore follows that $\rk_{O_{K,\fp^*}}(\ch{\ms}_{\fp^*}(K,T^*))
= 1$ as claimed.
\end{proof}

\begin{corollary} \label{C:vanishing}
Assume that $[\,,\,]_{K,\fp^*}$ is non-degenerate.

(a) If $r \geq 1$ and $\sha(K)(\fp^*)$ is finite, then
$$
\ord_{s=1} L_\fp^*(s)= r-1.
$$

(b) If $r=0$, then
$$
\ord_{s=1} L_\fp^*(s)=1.
$$
\end{corollary}

\begin{proof}
This follows directly from Propositions \ref{P:rkms} and
\ref{P:rkms2}, and \eqref{E:orders}.
\end{proof}

\begin{remark} Corollary \ref{C:vanishing}(b) confirms the
expectation expressed in \cite[Remark on p.74]{R1} (see also
\cite[\S11, Remarks(2)]{R}). It would be
interesting to know if there is any way of showing that
$\rk_{O_{K,\fp^*}}(\ms_{\fp^*}(K,T^*)) = 1$ when $r=0$
\textit{without} appealing to the functional equation satisfied by
$\cL_{\fp}$. \qed
\end{remark}

\begin{proposition} \label{P:sharel}
(a) Suppose that $r \geq 1$, and assume that $\sha(K)(\fp^*)$ is
finite. Then $\sha_{\rel(\fp)}(K)(\fp^*)$ is also finite, and we have
\begin{equation*}
|\sha_{\rel(\fp)}(K)(\fp^*)|
= |\sha(K)(\fp)| \cdot [ E(K_{\fp}) \otimes O_{K,\fp} :
\loc_\fp(\Sel(K,T))].
\end{equation*}

(b) Suppose that $r=0$. Then $\sha_{\rel(\fp)}(K)(\fp^*)$ has
$O_{K,\fp^*}$-corank one.
\end{proposition}

\begin{proof}
(a) For each $n \geq 1$, we define $B_n$ via exactness of
the sequence
$$
0 \to \sha(K)_{\pi^{*n}} \to H^1(K,E)_{\pi^{*n}} \to \prod_{v}
H^1(K_v, E)_{\pi^{*n}} \to B_n \to 0.
$$
Then there exists a map $h_n: H^1(K_{\fp},E)_{\pi^{*n}} \to B_n$,
and the sequence
\begin{equation} \label{E:seq1}
0 \to \sha(K)_{\pi^{*n}} \to \sha_{\rel(\fp)}(K)_{\pi^{*n}} \to
H^1(K_{\fp},E)_{\pi^{*n}} \xrightarrow{h_n} B_n
\end{equation}
is exact. Passing to direct limits over $n$ in \eqref{E:seq1} yields
the sequence
\begin{equation} \label{E:seq2}
0 \to \sha(K)(\fp^*) \to \sha_{\rel(\fp)}(K)(\fp^*) \to
H^1(K_\fp,E)(\fp^*) \xrightarrow{\varinjlim h_n} \varinjlim B_n.
\end{equation}
It follows from a theorem of Cassels (see \cite[p.198]{Ca}) that the
dual of $B_n$ is isomorphic to $\Sel(K,E_{\pi^n})$. Tate local duality
implies that the dual of $H^1(K_{\fp},E)_{\pi^{*n}}$ is isomorphic to
$E(K_{\fp})/\pi^n E(K_{\fp})$ and that the kernel of $\varinjlim h_n$
is isomorphic to the dual of the cokernel of the localisation map
$$
\loc_\fp: \ch{\Sel}(K,T) \to E(K_{\fp}) \otimes O_{K,\fp}.
$$
If $r \geq 1$, then this cokernel is finite, and we therefore deduce
that
$$
[\sha_{\rel(\fp)}(K)(\fp^*): \sha(K)(\fp^*)] = [ E(K_{\fp}) \otimes
O_{K,\fp} : \loc_{\fp}(\ch{\Sel}(K,T))].
$$
Hence, we have
\begin{equation*}
|\sha_{\rel(\fp)}(K)(\fp^*)| 
= |\sha(K)(\fp^*)| \cdot [ E(K_{\fp}) \otimes O_{K,\fp} :
\loc_{\fp}(\ch{\Sel}(K,T))]
\end{equation*}
as claimed.

(b) If $r=0$, then $\ch{\Sel}(K,T)$ is trivial, because $\sha(K)$ is
known to be finite, and $E(K)(\fp)=0$. This implies that
$\Coker(\loc_\fp) = E(K_\fp) \otimes O_{K,\fp} $ is $O_{K,\fp}$-free
of rank one. It now follows from \eqref{E:seq2} that
$\sha_{\rel(\fp)}(K)(\fp^*)$ has $O_{K,\fp^*}$-corank one.
\end{proof}

\begin{proposition} \label{P:shams}
Suppose that $r \geq 1$, and assume that $\sha(K)(\fp^*)$ is finite. Then
$$
|\ms_{\fp^*}(K,W^*)_{/\Div}| = |\sha_{\rel(\fp)}(K)(\fp^*)| \cdot
[ E(K_{\fp^*}) \otimes_{O_K} O_{K,\fp^*} :
\loc_{\fp^*}(\ch{\Sel}(K,T^*))].
$$
\end{proposition}

\begin{proof} 
Let $y_1,\ldots,y_{r-1}$ be an $O_{K,\fp^*}$-basis of
$E_{1,\fp^*}(K)$, and extend it to an $O_{K,\fp^*}$-basis
$y_1,\ldots,y_{r-1},y_{\fp^*}$ of $E(K) \otimes_{O_K}
O_{K,\fp^*}$. There is an exact sequence
$$
0 \to O_{K,\fp^*} \cdot y_{\fp^*} \to E(K_{\fp^*}) \otimes_{O_K}
O_{K,\fp^*} \to U \to 0,
$$
with
\begin{align*}
|U| &= [ E(K_{\fp^*}) \otimes_{O_K} O_{K,\fp^*}: \loc_{\fp^*}(E(K)
 \otimes_{O_K} O_{K,\fp^*})] \\
&= [ E(K_{\fp^*}) \otimes_{O_K} O_{K,\fp^*}:
 \loc_{\fp^*}(\ch{\Sel}(K,T^*))].
\end{align*}
Tensoring this sequence with $D_{\fp^*}$ yields an exact sequence
$$
0 \to V \to (O_{K,\fp^*} \cdot y_{\fp^*}) \otimes_{O_K} D_{\fp^*} \to
E(K_{\fp^*}) \otimes_{O_K} D_{\fp^*} \to 0,
$$
with $|U| = |V|$. As 
$$
E(K) \otimes _{O_K} O_{K,\fp^*} \simeq	E_{1,\fp^*}(K) \oplus
(O_{K,\fp^*} \cdot y_{\fp^*}),
$$
it follows that the kernel of the localisation map
$$
E(K) \otimes_{O_K} D_{\fp^*} \to E(K_{\fp^*}) \otimes_{O_K} D_{\fp^*}
$$
is isomorphic to $(E_{1,\fp^*}(K) \otimes_{O_K} D_{\fp^*}) \oplus V$.

Define
$$
\sha(K)_{\rel}:= \Ker \left[ H^1(K,E) \to \prod_{v \nmid p}
  H^1(K_v,E) \right];
$$
then we have an exact sequence
$$
0 \to E(K) \otimes D_{\fp^*} \to \Sel_{\rel}(K,W^*) \to
\sha_{\rel}(K)(\fp^*) \to 0.
$$
Now consider the following commutative diagram, in which the
vertical arrows are the obvious localisation maps:
\begin{equation*}
\begin{CD}
0 @>>> E(K) \otimes D_{\fp^*} @>>> \Sel_{\rel}(K,W^*) @>>>
\sha_{\rel}(K)(\fp^*) @>>> 0 \\
@. @VVV @VVV @VVV @. \\
0 @>>> E(K_{\fp^*}) \otimes D_{\fp^*} @>>>  H^1(K_{\fp^*},W^*) @>>>
H^1(K_{\fp^*},E)(\fp^*) @>>> 0
\end{CD}
\end{equation*}
Applying the Snake Lemma to this diagram yields the exact sequence
$$
0 \to (E_{1,\fp^*}(K) \otimes D_{\fp^*}) \oplus V \to \ms_{\fp^*}(K,W^*) \to
\sha_{\rel(\fp)}(K)(\fp^*) \to 0.
$$
As $\sha_{\rel}(K)(\fp^*)$ is finite (see Proposition \ref{P:sharel})
and $E_{1,\fp^*}(K) \otimes_{O_K} D_{\fp^*}$ is divisible, it follows
that
\begin{align*}
\ms_{\fp^*}(K,W^*)_{/\Div} &= |\sha_{\rel}(K)(\fp^*)| \cdot |V| \\
&=|\sha_{\rel(\fp)}(K)(\fp^*)| \cdot
[ E(K_{\fp^*}) \otimes_{O_K} O_{K,\fp^*} :
\loc_{\fp^*}(\ch{\Sel}(K,T^*))],
\end{align*}
as asserted.
\end{proof}

\section{Proof of Theorem \ref{T:A}} \label{S:ThmA}

\begin{proposition} \label{P:shams2}
Suppose that $r=0$. Then
$$
| \ms_{\fp^*}(K,W^*)_{/\Div}| \sim (1 - \psi(\fp^*)) \cdot
\frac{|\sha(K)_{\rel(\fp)}(\fp^*)_{/\Div}|}{[H^1(K_{\fp^*},T):
\loc_{\fp^*}(\ms_{\fp}(K,T))]}.
$$
\end{proposition}

\begin{proof}
Consider the following diagram in which all columns are exact and
$f_1$, $f_2$ are the obvious localisation maps:

\begin{equation*}
\begin{CD}
@.  0 @>>> \ms_{\fp^*}(K,W^*) @>>>
\sha_{\rel(\fp)}(K)(\fp^*) \\
@. @VVV   @VVV @VVV    \\
0 @>>> E(K) \otimes D_{\fp^*}=0 @>>> \Sel_{\rel}(K,W^*) @>>>
\sha_{\rel}(K)(\fp^*) @>>> 0 \\
@. @VVV @VV{f_1}V @VV{f_2}V @. \\
0 @>>> E(K_{\fp^*}) \otimes D_{\fp^*} @>>>  H^1(K_{\fp^*},W^*) @>>>
H^1(K_{\fp^*},E)(\fp^*) @>>> 0 \\
@. @VVV  @VVV  @VVV  @. \\
@. E(K_{\fp^*}) \otimes D_{\fp^*} @>>> \Coker(f_1) @>>> \Coker(f_2) @.
\end{CD}
\end{equation*}

Applying the Snake Lemma to this diagram yields an exact sequence
\begin{equation} \label{E:keyseq}
0 \to \ms_{\fp^*}(K,W^*) \to \sha_{\rel(\fp)}(K)(\fp^*) \to
E(K_{\fp^*}) \otimes D_{\fp^*} \to \Coker(f_1) \to \Coker(f_2) \to 0.
\end{equation}

Let us first determine $\Coker(f_1)$. The Poitou-Tate exact sequence
gives
$$
0 \to \ms_{\fp^*}(K,W^*) \to \Sel_{\rel}(K,W^*) \xrightarrow{f_1}
H^1(K_{\fp^*},W^*) 
\to \ch{\ms}_{\fp}(K,T)^\land \to H^2(G_{K,S_K},W^*),
$$
where $G_{K,S_K}$ denotes the Galois group over $K$ of the maximal
extension of $K$ that is unramified away from $p$. Since $r=0$, 
Propositions \ref{P:costr} and
\ref{P:rkstr1} imply that $H^2(G_{K,S_K},W^*) = 0$, and so we have
\begin{equation} \label{E:f1}
\Coker(f_1) \simeq \ch{\ms}_{\fp}(K,T)^\land.
\end{equation}
In particular, it follows from Lemma \ref{L:free} and Proposition
\ref{P:rkms2} that $\Coker(f_1)$ is divisible of $O_{K,\fp^*}$-corank
one.

In order to determine $\Coker(f_2)$, we observe that $E(K_{\fp^*})
\otimes D_{\fp^*}$ is divisible of $O_{K,\fp^*}$-corank one, and the
kernel of the map
$$
E(K_{\fp^*}) \otimes D_{\fp^*} \to \Coker(f_1)
$$
in \eqref{E:keyseq} is isomorphic to $\sha_{\rel(\fp)}(K)(\fp^*) /
\ms_{\fp^*}(K,W^*)$. This last group is finite, because both
$\sha_{\rel(\fp)}(K)(\fp^*)$ and $\ms_{\fp^*}(K,W^*)$ have
$O_{K,\fp^*}$-corank one (see Propositions \ref{P:sharel}(b) and
\ref{P:rkms2}). It therefore follows that $\Coker(f_2)=0$.

From \eqref{E:keyseq} and \eqref{E:f1}, we obtain the
sequence
\begin{equation} \label{E:keyseq2}
0 \to \frac{\sha_{\rel(\fp)}(K)(\fp^*)}{\ms_{\fp^*}(K,W^*)} \to
E(K_{\fp^*}) \otimes D_{\fp^*} \to \ch{\ms}_{\fp}(K,T)^\land \to 0.
\end{equation}
Dualising this sequence yields
$$
0 \to \ch{\ms}_{\fp}(K,T) \to
\frac{H^1(K_{\fp^*},T)}{H_f^1(K_{\fp^*},T)} \to
\left[\frac{\sha_{\rel(\fp)}(K)(\fp^*)}{\ms_{\fp^*}(K,W^*)}
\right]^\land \to 0.
$$
We therefore have
\begin{align*}
\left| \left[\frac{\sha_{\rel(\fp)}(K)(\fp^*)}{\ms_{\fp^*}(K,W^*)}
\right]^\land \right| &=
\left| \frac{\sha_{\rel(\fp)}(K)(\fp^*)}{\ms_{\fp^*}(K,W^*)} \right| \\
&=\left|
\frac{\sha_{\rel(\fp)}(K)(\fp^*)_{/\Div}}{\ms_{\fp^*}(K,W^*)_{/\Div}}
\right| \\
&= [ H^1(K_{\fp^*},T): \loc_{\fp^*}(\ch{\ms}_{\fp}(K,T))] \cdot |
H^1_f(K_{\fp^*},T)|^{-1},
\end{align*}
which in turn implies that
$$
| \ms_{\fp^*}(K,W^*)_{/\Div} | = \frac{ |
\sha_{\rel(\fp)}(K)(\fp^*)_{/\Div} |}{ [ H^1(K_{\fp^*},T):
\loc_{\fp^*}(\ch{\ms}_{\fp}(K,T))] } \cdot | H^1_f(K_{\fp^*},T)|.
$$
Since 
$$
| H^1_f(K_{\fp^*},T) | \sim 1 - \psi(\fp^*)
$$
(see Lemma \ref{L:Euler}), we finally obtain
$$
| \ms_{\fp^*}(K,W^*)_{/\Div}| \sim
(1 - \psi(\fp^*)) \cdot
\frac{|\sha(K)_{\rel(\fp)}(\fp^*)_{/\Div}|}{[H^1(K_{\fp^*},T):
\loc_{\fp^*}(\ms_{\fp}(K,T))]},
$$
as claimed.
\end{proof}

{\it Proof of Theorem \ref{T:A}.} We first note that, as
$[\,,\,]_{K,\fp^*}$ is non-degenerate (by hypothesis), we have
$\ord_{s=1} L_{\fp}^{*}(s) = 1$ (Corollary
\ref{C:vanishing}(b)). Hence from \eqref{E:leading}, \eqref{E:LT6},
Proposition \ref{P:shams2} and Remark \ref{R:char}, we have
\begin{align*} 
\lim_{s \to 1} \frac{L_\fp^*(s)}{s-1} &\sim \log_p(\psi^*(\gamma))
\cdot \frac{H_K}{t} \Biggr |_{t=0}  \\
&\sim \log_p(\psi^*(\gamma)) \cdot \left|
\ms_{\fp^*}(K,W^*)_{/\Div} \right| \cdot
\cR_{K,\fp^*} \\
&\sim \log_p(\psi^*(\gamma)) \cdot
(1 - \psi(\fp^*) \cdot
\frac{|\sha_{\rel(\fp)}(K)(\fp^*)_{/\Div}|}{[H^1(K_{\fp^*},T):
\loc_{\fp^*}(\ms_{\fp}(K,T))]} \cdot \cR_{K,\fp^*}.
\end{align*}
This completes the proof of Theorem \ref{T:A}. \qed

\section{Proof of Theorem \ref{T:B}} \label{S:ThmB}

Suppose now that $r \geq 1$.  Then $E(K) \otimes O_{K,\fp^*}$ is a
free $O_{K,\fp^*}$-module of rank $r$. Proposition \ref{P:rkstr1}
implies that the kernel of the localisation map
$$
\loc_{\fp^*}: E(K) \otimes_{O_K} O_{K,\fp^*} \to E(K_{\fp^*}) \otimes
O_{K,\fp^*}
$$
has $O_{K,\fp^*}$-rank $r-1$. Let $y_1,\ldots,y_{r-1}$ be an
$O_{K,\fp^*}$-basis of this kernel, and extend it to an
$O_{K,\fp^*}$-basis $y_1,\ldots,y_{r-1},y_{\fp^*}$ of $E(K) \otimes
O_{K,\fp^*}$. 

\begin{proposition} \label{P:index}
With the above assumptions and notation, we have
$$
[ E(K_{\fp^*}) \otimes_{O_K} O_{K,\fp^*}: \loc_{\fp^*}(E(K)
\otimes_{O_K} O_{K,\fp^*})] 
\sim
p^{-1} \log_{E,\fp^*}(y_{\fp^*}),
$$
where $\log_{E,\fp^*}$ denotes the $\fp^*$-adic logarithm associated to
$E$. Similarly, we also have
$$
[ E(K_{\fp}) \otimes_{O_K} O_{K,\fp}: \loc_{\fp}(E(K) \otimes_{O_K}
O_{K,\fp})] \sim p^{-1} \log_{E,\fp}(y_{\fp}),
$$
when $y_{\fp} \in E(K_{\fp}) \otimes_{O_K} O_{K,\fp}$ is defined
analogously to $y_{\fp^*}$.
\end{proposition}

\begin{proof}
We give the proof of the first assertion; that of the second is of
course essentially identical.

We first observe that, from the definitions, we have
$$
[ E(K_{\fp^*}) \otimes_{O_K} O_{K,\fp^*}: \loc_{\fp^*}(E(K)
\otimes_{O_K} O_{K,\fp^*})]=
[ E(K_{\fp^*}) \otimes O_{K,\fp^*}: \loc_{\fp^*}(O_{K,\fp^*} \cdot y_{\fp^*})].
$$

Let $E_0$ denote the kernel of reduction modulo $\fp^*$ of $E$, so we
have an exact sequence
$$
0 \to E_0(K_{\fp^*}) \to E(K_{\fp^*}) \to \tilde{E}_{\fp^*}(k_{\fp^*}) \to 0.
$$
Set
$$
Z:= O_{K,\fp^*} \cdot y_{\fp^*},\quad Z_0:= \loc_{\fp^*}(Z) \cap
E_0(K_{\fp^*}), \quad C:= \loc_{\fp^*}(Z)/Z_0.
$$
Write $\lambda_{\fp^*}$ for the restriction of $\loc_{\fp^*}$ to $Z$. We
have the following commutative diagram:
\begin{equation*}
\begin{CD}
0 @>>> Z_0  @>>>  Z  @>>> C
\otimes_{O_K} O_{K,\fp^*} @>>> 0 \\
@. @VV{\rho}V  @VV{\lambda_{\fp^*}}V  @VV{\rho'}V  @. \\
0 @>>> E_0(K_\fp^*) \otimes_{O_K} O_{K,\fp^*} @>>> E(K_{\fp^*})
\otimes_{O_K} O_{K,\fp^*} 
@>>> \tilde{E}_{\fp^*}(k_{\fp^*}) 
\otimes_{O_K} O_{K,\fp^*} @>>> 0
\end{CD}
\end{equation*}

Observe that $\rho$ is injective since $\lambda_{\fp^*}$ is injective,
and that $\tilde{E}_{\fp^*}(k_{\fp^*}) \otimes_{O_K} O_{K,\fp^*} = 0$
because $\tilde{E}_{\fp^*}(k_{\fp^*})(p)
=\tilde{E}_{\fp^*}(k_{\fp^*})(\fp)$  (see e.g. \cite[p. 28]{PR7}). 
Applying the Snake Lemma to the diagram yields the exact
sequence
$$
0 \to \Ker(\rho') \to \Coker(\rho) \to \Coker(\lambda_{\fp^*}) \to 0,
$$
and so we have
$$
|\Coker(\lambda_{\fp^*})| = |C \otimes_{O_K} O_{K,\fp^*}|^{-1} \cdot |
\Coker(\rho) |.
$$

Set $k = [Z:Z_0] = |C \otimes O_{K,\fp^*}|$; then $ky_{\fp^*}$ is
an $O_{K,\fp^*}$-generator of $Z_0$. Since there is an isomorphism
$$
\log_{E,\fp^*}: E_0(K_{\fp^*}) \xrightarrow{\sim}
\fp^* O_{K,\fp^*},
$$
it follows that we have
$$
| \Coker(\rho) | \sim p^{-1} \log_{E,\fp^*}(ky_{\fp^*}) = k p^{-1}
  \log_{E,\fp^*}(y_{\fp^*}).
$$
Therefore
$$
|\Coker(\lambda_{\fp^*})| \sim  p^{-1} \log_{E,\fp^*}(y_{\fp^*}),
$$
and this establishes the desired result.
\end{proof}

\begin{corollary} \label{C:ordsharel}
Suppose that $r \geq 1$ and assume that $\sha(K)(\fp^*)$ is
finite. Then
\begin{equation*}
|\sha_{\rel(\fp)}(K)(\fp^*)| 
= p^{-1}\cdot |\sha(K)(\fp^*)| \cdot \log_{E,\fp}(y_{\fp}).
\end{equation*}
\end{corollary}

\begin{proof} This follows directly from Propositions
\ref{P:sharel}(a) and \ref{P:index}.
\end{proof}

{\it Proof of Theorem \ref{T:B}.} By hypothesis, $[\,
,\,]_{K,\fp^*}$ is non-degenerate, $r \geq 1$, and $\sha(K)(p)$ is
finite; hence we have that $\ord_{s=1} L_{\fp}^{*}(s) = r-1$
(Corollary \ref{C:vanishing}(a)). Proposition \ref{P:shams} and
Corollary \ref{C:ordsharel} imply that
\begin{align*}
| \ms_{\fp^*}(K,W^*)_{/\Div}| &= |\sha_{\rel(\fp)}(K)(\fp^*)| \cdot
[E(K_{\fp^*}) \otimes_{O_K} O_{K,\fp^*}: \loc_{\fp^*}
(\ch{\Sel}(K,T^*))] \\
&\sim p^{-2}\cdot |\sha(K)(\fp^*)| 
\cdot \log_{E,\fp^*}(y_{\fp^*}) \cdot \log_{E,\fp}(y_{\fp}) .
\end{align*}
We therefore deduce from \eqref{E:leading}, \eqref{E:LT6} and Remark
\ref{R:char} that
\begin{align*}
\lim_{s \to 1} &\frac{L_\fp^*(s)}{(s-1)^{r-1}} \sim \\
& [\log_p(\psi^*(\gamma))]^{r-1} \cdot
p^{-2} \cdot
|\sha(K)(\fp^*)| \cdot
\log_{E,\fp^*}(y_{\fp^*}) \cdot \log_{E,\fp}(y_{\fp})
\cdot \cR_{K,\fp^*},
\end{align*}
as asserted.

This completes the proof of Theorem \ref{T:B}. \qed

\section{Canonical elements in restricted Selmer groups}
\label{S:exact}

The goal of this section is to explain how the methods of \cite{R} may
be used to produce an exact formula for $\lim_{s \to 1}
L_{\fp}^{*}(s)/(s-1)$ when $r=0$ (see Theorem \ref{T:exact} below). The
arguments involved are quite similar to those of \cite{R}, and so, in
what follows, we assume that the reader has a copy of \cite{R} and is
willing to refer to it from time to time for some of the details we omit.

We begin by introducing the following notation (some of which differs
from that of \cite{R}):
\begin{align*}
&U_{n,\fp}:= \text{units in $\cK_{n,\fp}$ congruent to $1$ modulo
    $\fp$;} \\
&U_{n,\fp^*}:= \text{units in $\cK_{n,\fp^*}$ congruent to $1$ modulo
    $\fp^*$;} \\
&U_{\infty,\fp}:= \varprojlim U_{n,\fp},\quad U_{\infty,\fp^*}:= 
\varprojlim U_{n,\fp^*};\\
&U^{*}_{n,\fp}:= \text{units in $\cK^{*}_{n,\fp}$ congruent to $1$ modulo
    $\fp$;} \\
&U^{*}_{n,\fp^*}:= \text{units in $\cK^{*}_{n,\fp^*}$ congruent to $1$ modulo
    $\fp^*$;} \\
&U^{*}_{\infty,\fp}:= \varprojlim U^{*}_{n,\fp},\quad U^{*}_{\infty,\fp^*}:=
    \varprojlim U_{n,\fp^*},
\end{align*}
where all inverse limits are taken with respect to norm maps. We also set
\begin{align*}
&\cE_n:= \text{global units of $\cK_n$},\quad \cE_n^*:= \text{global
  units of $\cK_n^*$};\\
&\ov{\cE}_n:= \text{the closure of the projection of $\cE_n$ into
  $U_{n,\fp}$}; \\
&\ov{\cE}^*_n:= \text{the closure of the projection of $\cE^*_n$ into
  $U^*_{n,\fp^*}$}; \\
&\ov{\cE}_{\infty}:= \varprojlim \ov{\cE}_n, \quad \ov{\cE}^{*}_{\infty}:=
  \varprojlim \ov{\cE}^*_n.
\end{align*}

\begin{remark} \label{R:leo}
Note that since the strong Leopoldt conjecture holds for all
abelian extensions of $K$ (see \cite{B}), we have that
$$
\ov{\cE}_n \simeq \ov{\cE}_n \otimes_{\bZ} \bZ_p, \quad 
\ov{\cE}^*_n \simeq \ov{\cE}^*_n \otimes_{\bZ} \bZ_p,
$$ and so we may also view $\ov{\cE}_{\infty}$ as being a submodule of
$U_{\infty,\fp^*}$ and $\ov{\cE}^{*}_{\infty}$ as being a submodule of
$U^{*}_{\infty,\fp}$. We shall do this without further comment several
times in what follows.  \qed
\end{remark}

\begin{proposition} \label{P:kummerinj}
There are natural injections
\begin{align*}
&\rho: \Hom(T^*,(U^{*}_{\infty,\fp} \otimes
\bQ)/\ov{\cE}^*_\infty)^{\Gal(\cK_\infty^*/K)} \hookrightarrow
\ch{\ms}_{\fp}(K,T), \\
&\rho^*: \Hom(T,(U_{\infty,\fp} \otimes
\bQ)/\ov{\cE}_\infty)^{\Gal(\cK_\infty/K)} \hookrightarrow
\ch{\ms}_{\fp^*}(K,T^*)
\end{align*}
\end{proposition}

\begin{proof}
The proof of this result is essentially the same, \textit{mutatis
  mutandis}, as that of \cite[Proposition 2.4]{R}. The map $\rho$ is
  defined as follows.

For any $f \in \Hom(T^*,(U^{*}_{\infty,\fp} \otimes
\bQ)/\ov{\cE}^*_\infty)^{\Gal(\cK_\infty^*/K)}$ and any integer $n
\geq 1$, we define $f_n \in \Hom(E_{\pi^n},
\cE^{*}_{n}/\cE_{n}^{*p^n})^{\Gal(\cK_\infty/K)}$ to be the image of $f$
under the following composition of maps:
\begin{align*}
\Hom(T^*,(U^{*}_{\infty,\fp} \otimes
\bQ)/\ov{\cE}^*_\infty)^{\Gal(\cK_\infty^*/K)} &\to 
\Hom(T^*,(U^{*}_{n,\fp} \otimes
\bQ)/\ov{\cE}^*_n)^{\Gal(\cK_\infty^*/K)} \\
&\to \Hom(E_{\pi^n}, \cE^*_n/\cE_{n}^{*p^n})^{\Gal(\cK^*_\infty/K)},
\end{align*}
where the first arrow is the map induced by the natural projection
$U^{*}_{\infty,\fp} \to U^{*}_{n,\fp}$, and the second arrow is
induced by raising to the $p^n$-th power in $U^{*}_{n,\fp}$.

Recall that, for each $n \geq 1$, there is an isomorphism
$$
\rho_n: H^1(K,E_{\pi^n}) \xrightarrow{\sim} \Hom(E_{\pi^{*n}},
\cK_{n}^{*\times}/\cK_{n}^{*\times p^{n}})^{\Gal(\cK_n^*/K)}
$$
(see e.g. \cite[Lemma 2.1]{R} or \cite[Lemme 12]{PR1}). We define
$$
\rho(f):= [(p-1)(\pi^{*})^{2n} \rho_{n}^{-1}(f_n)] \in \varprojlim_n
H^1(K,E_{\pi^n}).
$$
It is not hard to check from the definition that $\rho$ is
injective. It follows from Theorem \ref{T:coates}, Proposition
\ref{P:control}, and Corollary \ref{C:descent} that
$\rho_{n}^{-1}(f_n) \in \ms_{\fp}(K,E_{\pi^n})$ if and only if the
restriction of $\rho_{n}^{-1}(f_n)$ to $H^1(\fK_\infty,E_{\pi^n})$ is
unramified outside $\fp^*$. It may be shown via an argument very
similar to that given in \cite[Lemmas 2.1 and 2.3]{R} that this in
fact the case. 
\end{proof}

We shall now explain how elliptic units may be used (following
\cite{R}) to construct canonical elements
$$
s_{\fp}^{(1)} \in \ch{\ms}_{\fp}(K,T),\quad s_{\fp^*}^{(1)} \in
\ch{\ms}_{\fp^*}(K,T^*)
$$
when $r=0$. These are the analogues in the present situation of the
elements $x_{\fp}^{(1)} \in \ch{\Sel}(K,T)$ and $x_{\fp^*}^{(1)} \in
\ch{\Sel}(K,T^*)$ constructed in \cite{R} when $r=1$.

Let $\cC_\infty \subseteq \cE_\infty$ and $\cC^*_\infty \subseteq
\cE^*_\infty$ denote the norm-coherent systems of elliptic units
constructed in \cite[\S3]{R}, and write $\ov{\cC}_\infty$
and $\ov{\cC}^*_\infty$ for the closure of $\cC_\infty$ in
$\ov{\cE}_\infty$ and $\cC^*_\infty$ in $\cE^*_\infty$
respectively. Set
$$
\cJ^*:= \Ker(\psi^*: \Lambda(\cK^*_\infty) \to \bZ_p), \quad
\cJ:= \Ker(\psi: \Lambda(\cK_\infty) \to \bZ_p),
$$
and let $\vt^*$ be the generator of $\cJ^*$ fixed in
\cite[\S6]{R} (so $\vt^* = \gamma \psi^*(\gamma^{-1}) -1$, where
$\gamma$ is any topological generator of $\Gal(\cK_\infty^*/K)$
satisfying $\log_p(\psi^*(\gamma)) = p$). Write $\ff \subseteq O_K$
for the conductor of the Grossencharacter associated to $E$, and let
$\bN(\ff)$ denote the norm of this ideal. Fix $B \in
E_{\ff}/\Gal(\ov{K}/K)$, and generators $w$ of $T$ and $w^*$ of $T^*$
according to the recipe described in \cite[\S6]{R}. Let
$$
\theta_B(\bN(\ff)^{-1}w^*) \in \ov{\cC}^*_\infty \subseteq
U_{\infty,\fp}^{*} \otimes \bQ
$$ denote the elliptic unit constructed in \cite[\S3]{R}.

Suppose that $t$ is a positive integer such that
\begin{equation*}
\ov{\cC}^*_\infty \subseteq \cI^{t-1} \ov{\cE}_\infty^* \subseteq
U_{\infty,\fp}^{*} \otimes \bQ \quad \text{and} \quad
\ov{\cC}^*_\infty \subseteq \cI^{t}(U_{\infty,\fp}^{*} \otimes \bQ).
\end{equation*}

\begin{proposition} \label{P:canhom}
There exists a unique homomorphism $\sigma_{\fp}^{(t)} \in
\Hom(T^*,(U_{\infty,\fp}^{*} \otimes \bQ)/\ov{\cE}_\infty^*)$ such that
$$
\sigma_{\fp}^{(t)}(w^*)^{\vt^{*t}} = \theta_B(-\bN(\ff)^{-1}w^*)
$$
in $\ov{\cE}_\infty^*/\cJ^{*t}\ov{\cE}_\infty^*$.
\end{proposition}

\begin{proof} Theorem 7.2(i) of \cite{R} implies that
  $U_{\infty,\fp}^{*}$ contains no $\vt^*$-torsion elements. The
  existence of $\sigma_{\fp}^{(t)}$ therefore follows via an argument
  very similar to that of \cite[Theorem 4.2]{R}.
\end{proof}

We set
$$
s_{\fp}^{(t)}:= \rho(\sigma_{\fp}^{(t)}),\quad s_{\fp^*}^{(t)}:=
\rho^*(\sigma_{\fp^*}^{(t)}),
$$ 
where of course the definition $\sigma_{\fp^*}^{(t)} \in
\Hom(T,(U_{\infty,\fp^*} \otimes \bQ)/\ov{\cE}_{\infty})$ the same,
\textit{mutatis mutandis}, as that of $\sigma_{\fp}^{(t)}$.

\begin{remark} In fact the only non-zero values of $s_{\fp}^{(t)}$ and
  $s_{\fp^*}^{(t)}$ occur when $r=0$ and $t=1$:

(a) Suppose that $r=0$. Then $L_{\fp}(1) \neq 0$, and so we have (via
  \cite[Theorem 7.2(i)]{R}, for example):
$$
\ov{\cC}_\infty \subseteq \ov{\cE}_\infty \subset U_{\infty,\fp}
\otimes \bQ \quad \text{and} \quad \ov{\cC}_\infty \not \subseteq
\cI(U_{\infty,\fp} \otimes \bQ).
$$
In particular, we have that $\ov{\cC}_\infty \not \subseteq \cI
\ov{\cE}_\infty \subseteq U_{\infty,\fp} \otimes \bQ$. Similar remarks
imply that also $\ov{\cC}^*_\infty \not \subseteq \cI^*
\ov{\cE}^*_\infty \subseteq U^{*}_{\infty,\fp^*} \otimes
\bQ$. Applying Remark \ref{R:leo}, we deduce that
\begin{equation} \label{E:r1}
\ov{\cC}^*_\infty \not \subseteq \cI^* \ov{\cE}^*_\infty \subseteq
U_{\infty,\fp^*}^{*} \otimes \bQ.
\end{equation} 
Now suppose in addition that $[\,,\,]_{K,\fp^*}$ is
non-degenerate. Then Theorem \ref{T:A} implies that $\ord_{s=1}
L_{\fp}^{*}(s) = 1$, and so from \cite[Theorem 7.2(i)]{R}, we have
\begin{equation} \label{E:r2}
\ov{\cC}_\infty^* \subseteq \cI^* (U_{\infty,\fp}^{*} \otimes \bQ).
\end{equation}
We now deduce from \eqref{E:r1} and \eqref{E:r2} and the definition of
$\rho$ that $s_{\fp}^{(1)} \neq 0$.

A similar argument shows that $s_{\fp^*}^{(1)} \neq 0$ also.
\smallskip

(b) Suppose now that $r \geq 1$. Assume that $\sha(K)(p)$ is finite,
and that the height pairing $[\,,\,]_{K,\fp^*}$ is
non-degenerate. Then Theorem \ref{T:B} (or \cite[Corollary 11.3]{R})
implies that $\ord_{s=1}
L_{\fp}^{*}(s) = r-1$, and so it follows from \cite[Theorem 7.2(i)]{R}
that
\begin{equation} \label{E:r3}
\ov{\cC}^*_\infty \subseteq \cI^{*r-1}(U_{\infty,\fp}^{*} \otimes
\bQ).
\end{equation}

On the other hand, Theorem 4.2 and Proposition 4.4 of \cite{R} imply
that 
\begin{equation*}
\ov{\cC}_\infty^* \subseteq \cI^{*r-1}\ov{\cE}^*_\infty \subseteq 
U_{\infty,\fp^*}^{*} \otimes \bQ, \quad
\ov{\cC}_\infty^* \not \subseteq \cI^{*r}\ov{\cE}^*_\infty \subseteq 
U_{\infty,\fp^*}^{*} \otimes \bQ,
\end{equation*}
and so applying Remark \ref{R:leo}, we deduce that
\begin{equation} \label{E:r4}
\ov{\cC}_\infty^* \subseteq \cI^{*r-1}\ov{\cE}^*_\infty \subseteq 
U_{\infty,\fp}^{*} \otimes \bQ, \quad
\ov{\cC}_\infty^* \not \subseteq \cI^{*r}\ov{\cE}^*_\infty \subseteq 
U_{\infty,\fp}^{*} \otimes \bQ.
\end{equation}

It now follows from \eqref{E:r3} and \eqref{E:r4} that $s_{\fp}^{(t)}
= 0$ for $1 \leq t \leq r-2$ and that $s_{\fp}^{(t)}$ is not defined
for $t \geq r-1$.
\smallskip

(c) Suppose that $r=0$, but that $\ord_{s=1} L_{\fp}^{*}(s)>1$ (so, in
particular, the pairing $[\,,\,]_{K,\fp^*}$ is degenerate, which we
expect never to happen). Then an argument similar to that given in (b)
above shows that $s_{\fp}^{(1)} = 0$, and that $s_{\fp}^{(t)}$ is not
defined for $t>1$. \qed
\end{remark}

\begin{theorem} \label{T:exact}
Suppose that $r=0$ and that $[\,,\,]_{K,\fp^*}$ is non-degenerate, so
$\ord_{s=1}L_{\fp}^{*}(s) = 1$. Then
$$
\lim_{s \to 1} \frac{L_{\fp}^{*}(s)}{s-1} =
\bN(\ff)^{-1}  (p-1) \left( 1 - \frac{\psi^*(\fp)}{p} \right) \lim_{n
  \to \infty} \log_{\fp}(\sigma_{\fp,n}^{(1)}(w^*)).
$$
\end{theorem}

\begin{proof} This may be shown in exactly the same way as
  \cite[Proposition 9.4(ii)]{R}.
\end{proof} 

\begin{remark} The precise relationship between Theorem \ref{T:A} and
  Theorem \ref{T:exact} is not clear, and it would be interesting to
  obtain a better understanding of this. 
\qed
\end{remark}

\end{document}